\documentclass[a4paper]{article}

\title{
A set optimization approach to zero-sum matrix games with multi-dimensional payoffs\footnote{The work of the first author was supported by a generous start-up grant from Free University Bozen. The authors are grateful for comments, hints and suggestions made by L. Andreozzi and the participants of the Economics Research seminar at Free University Bozen--in particular F. Boffa who initiated the seminar discussion--as well as by T. Tanaka.}
}

\author{
Andreas H. Hamel\footnote{Free University Bozen-Bolzano, Italy, Faculty of Economics and Management, \href{mailto:andreas.hamel@unibz.it}{andreas.hamel@unibz.it}}, Andreas L\"ohne\footnote{Friedrich Schiller University Jena, Germany, Department of Mathematics, \href{mailto:andreas.loehne@uni-jena.de}{andreas.loehne@uni-jena.de}}
}
\date{{\small \today}}

\usepackage{array} 
\usepackage{verbatim} 

\usepackage{amsmath, amssymb, enumerate, color, a4, booktabs, graphicx, enumerate}

\newtheorem{theorem}{Theorem}

\newtheorem{remark}[theorem]{Remark}
\newtheorem{lemma}[theorem]{Lemma}
\newtheorem{definition}[theorem]{Definition}
\newtheorem{proposition}[theorem]{Proposition}
\newtheorem{example}[theorem]{Example}

\numberwithin{equation}{section}  
\numberwithin{figure}{section}    
\numberwithin{table}{section}     
\numberwithin{theorem}{section}

\newcommand{\of}[1]{\ensuremath{\left( #1 \right)}}

\newcommand{\cb}[1]{\ensuremath{ \left\{ #1 \right\} }}
\newcommand{\sqb}[1]{\ensuremath{ \left[ #1 \right] }}

\newcommand{\bs}{\backslash}

\newcommand{\pend}{ \hfill $\square$ \medskip}

\renewcommand{\P}{\ensuremath{\mathcal{P}}}

\newcommand{\G}{\ensuremath{\mathcal{G}}}

\newcommand{\lel}{\preccurlyeq}
\newcommand{\leu}{\curlyeqprec}
\newcommand{\Min}{{\rm Min\,}}
\newcommand{\Max}{{\rm Max\,}}
\newcommand{\WMIN}{{\rm wMin\,}}
\newcommand{\WMAX}{{\rm wMax\,}}
\newcommand{\MIN}{{\text{\rm MIN$(I)$}}}
\newcommand{\MAX}{{\text{\rm MAX$(II)$}}}

\newcommand{\V}{\mathcal{V}}
\newcommand{\W}{\mathcal{W}}
\newcommand{\T}{\mathcal{T}}

\newcommand{\R}{\mathbb{R}}
\newcommand{\Z}{\mathbb{Z}}

\newcommand{\cl}{{\rm cl \,}}

\newcommand{\co}{{\rm co \,}}

\newcommand{\Int}{{\rm int\,}}

\usepackage[colorlinks=true, linkcolor=blue, citecolor=blue, urlcolor=blue, filecolor=blue]{hyperref}

\definecolor{color0}{gray}{.50}
\definecolor{color1}{rgb}{0,.2,.8}
\definecolor{color2}{rgb}{1,.2,0}
\definecolor{color3}{rgb}{.8,.5,1}

\begin{document}
\maketitle

\begin{abstract} A new solution concept for two-player zero-sum matrix games with multi-dimensional payoff is introduced. It is based on extensions of the vector order in $\R^K$ to order relations in the power set of $\R^K$, so-called set relations, and strictly motivated by the interpretation of the payoff as multi-dimensional loss for one and gain for the other player. The new concept provides coherent worst case estimates, i.e. minimax and maximin strategies, for games with multi-dimensional payoffs. It is shown that--in contrast to games with one-dimensional payoffs--minimax/maximin strategies are independent from Shapley's notion of equilibrium strategies for such games. Therefore, the two concepts are combined into new equilibrium notions for which existence theorems are given. By means of examples, relationships of the new concepts to existing ones such as Shapley and vector equilibria, vector minimax/maximin solutions and Pareto optimal security strategies are clarified. A algorithm for computing optimal strategies is presented.
\end{abstract}




\section{Introduction}
\label{SecIntro}

This note is an attempt to provide novel answers to the question what should actually be played in games in which the players have incomplete preferences for the payoffs, i.e. there are non-comparable alternatives. Our findings indicate that previously defined extensions of Nash equilibria to the case of vector-valued payoffs or utility functions are not the only natural path to follow as done in \cite[p. 311]{Bade05ET}: `There is a natural way of extending the standard notion of Nash equilibrium' to incomplete preferences. In particular, we will discuss that there might be other `incentives to deviate from her own action given every one else's action' (\cite[p. 311]{Bade05ET}) than just (strictly) improving the multi-dimensional payoff/utility with respect to the underlying vector order or (non-complete) preference.

Incomplete preferences and their representation came more into focus of the economic literature in the past two decades. Aumann \cite{Aumann62Economet} (in search for a utility representation theory without the completeness axiom as the `most questionable') and Bewley \cite{Bewley86R} (in search for a model that would make precise the difference between risk and uncertainty) raised the issue. Compare Ok et al \cite{OkOrtolevaRiella12Economet} for a synthesis of these approaches. Representation results for incomplete preferences via families of utility functions were given by Ok \cite{Ok02JET}, Dubra et al \cite{DubraMaccheroniOk04JET}, Evren, Ok \cite{EvrenOk11JME}, Bosi, Herden \cite{BosiHerden12JME} along with the perspective that vector optimization might be useful for solving the resulting multi-utility maximization problems: `Moreover, as it reduces finding the maximal elements in a given subset of $X$ with respect to $\preceq$ to a multi-objective optimization problem (cf. cf. Evren and Ok, 2011), in applications, this approach is likely to be more useful than the Richter-Peleg approach.' However, to the best of our knowledge, up to now there is no (expected) multi-utility maximization theory with results comparable to the scalar case. Rather, even in recent applications, a real-valued multi-attribute utility function is assumed to model the preference of a decision maker even in the highly complex framework of financial markets with frictions (cf. Campi, Owen \cite{CampiOwen11FS}).

A major motivation for investigating incomplete preferences stems from social choice theory. The reader may compare Eliaz and Ok \cite{EliazOk06GEB} (and the references therein) where the difference between indifference and indecisiveness is discussed as well as the weak axiom of revealed preferences (WARP) which almost forces a complete preference. It is shown that a weakened version of WARP is perfectly compatible with incomplete preferences, and the methodology given in this paper admits to distinguish between cases in which the decision maker is indifferent and cases in which (s)he is indecisive. 

The theory proposed here gives a different perspective to the very same phenomena. It is shown that players may deviate from an equilibrium by ``free will," i.e. they may want to move to another equilibrium which produces a payoff which is not comparable to the first for reasons which are not part of the mathematical model. Seen in this way, the interpretation is very close to a literal understanding of the word ``incomplete" as, for example, discussed in Carlier, Dana \cite{CarlierDana13JET}. Our novel contribution can be summarized as follows: ``Set relations" (extensions of the vector order for the payoffs to the power set of possible payoffs) are used for giving the players a guidance through the jungle of the many (non-comparable) equilibrium values with respect to the vector order, and, for the first time, well-motivated worst case estimates are provided which correspond to the von Neumann minimax and maximin strategies in the one-dimensional payoff case.

We use one of the easiest problems involving incomplete preferences as a show case for the set relation approach to game theory, namely a two-player, zero-sum matrix game where the utility function for each player just is, admitting mixed strategies, her/his expected vector-valued gain. 

Already Shapley \cite[p. 58]{Shapley59NRLQ} (see also Aumann \cite[p. 447]{Aumann62Economet}) gave a motivation for studying games with multi-dimensional payoffs: `The payoff of a game sometimes most naturally takes the form of a vector having numerical components that represent commodities (such as men, ships, money, etc.) whose relative values cannot be ascertained. The utility spaces of the players can therefore be given only a partial ordering (representable as the intersection of a finite number of total orderings), and the usual notions of solution must be generalized.' However, the theory even for ``simple" zero-sum matrix games with multi-dimensional payoffs is far from being complete in a way which parallels the elegant von Neumann approach in the scalar case. The earliest proposal by Blackwell \cite{Blackwell56PJM} generalized the notion of the value of a game to an ``approachable set" apparently motivated by the fact that, in contrast to the classic von Neumann situation, insisting on a single payoff (vector) as value of the game does not make much sense for games with vector payoffs.

Shapley \cite{Shapley59NRLQ} introduced an equilibrium notion which relies on the order generated by the positive cone in the images space. His concept indeed admitted to `utilize the theory of vector optimization (multi-objective programming)' (cf. Ok \cite[p. 433]{Ok02JET}), and it was put into the context of vector optimization as vector equilibrium points or vector saddle points, see Nieuwenhuis \cite[Definition 3.1]{Nieuwenhuis83JOTA}, Corley \cite[(5)]{Corley85JOTA}, Luc and Vargas \cite{LucVargas92NA}, Tanaka \cite{Tanaka94JOTA} as well as many others. Compare the survey Tanaka \cite{Tanaka00Incoll} and the references therein as well as Zhao \cite{Zhao91IJGT} where general multi-objective games are considered and more references can be found.

However, it turned out that, despite this effort, vector optimization did not yet provide a satisfying equilibrium theory, not even for two-player, zero-sum matrix games with vector-valued payoffs/utilities.   

What are the reasons for the lack of an appealing and applicable, multi-dimensional equilibrium theory under incomplete preferences? 

One difficulty is the sheer number of ``Shapley equilibrium values." Shapley's main result \cite[Theorem]{Shapley59NRLQ} reduces a zero-sum matrix game with vector payoffs to infinitely many scalar non-zero-sum games. In most cases, this produces just too many candidates for a solution and `the impossible task of solving all possible scalarizations' as remarked by Corley \cite[Section 4]{Corley85JOTA}. In a more general context, Bade \cite[p. 328]{Bade05ET} observed the same: `Since Nash equilibrium sets of games with incomplete preferences can be large, it is of interest to consider refinements of equilibria $\ldots$'  We refer to De Marco, Morgan \cite{DeMarcoMorgan07IGTR} for a discussion of several approaches and references. The idea is to impose additional requirements to the equilibrium \cite[page 171]{DeMarcoMorgan07IGTR} `which adds to the original problem new endogenous parameters that are typical for the vector-valued form.' However, the questions then asked are almost exclusively the one for the existence of such refined equilibria and the one for stability with respect to some class of perturbations. In the context of vector optimization, similar scalarization results can be found, for example, in Zeleny \cite{Zeleny75IJGT} (basically, already Shapley's result), Cook \cite{Cook76NRLQ} (minimizing the weighted underachievement of goals) and Wierzbicki \cite{Wierzbicki95JSEE} (via nonlinear scalarization functions). Zhao's notion of Nash equilibria \cite[Definition 12]{Zhao91IJGT} is based on properly efficient points which can also be found via linear scalarizations.

Another, even more important difficulty is that Shapley equilibria are missing two main features of equilibrium strategies for scalar games. They are not interchangeable, and they do not provide worst case estimates.

The missing interchangeability was observed by Corley \cite[Example 3.2]{Corley85JOTA} (see Example \ref{ExCorley3.2} below), and contradicts Aumann's belief that `the interchangeability property holds' (cf. Aumann \cite[p. 455]{Aumann62Economet}) for Shapley equilibrium strategies. Therefore, a ``best possible answer" to an equilibrium strategy is not necessarily a ``best possible answer" to another which leaves the question completely open what kind of strategy protects a player best against the opponent's choices. Below, the meaning of ``best" will be scrutinized, and it will be given a new meaning via a set relation approach.

The missing worst case insurance motivated (sightly deviating from each other) concepts of vector minimax and maximin strategies as discussed, for example, by Nieuwenhuis \cite{Nieuwenhuis83JOTA}, Corley \cite{Corley85JOTA} and Tanaka \cite{Tanaka88JOTA}. The relationship to the previously introduced vector equilibrium/saddle points remained obscure as observed by Corley \cite[Section 4]{Corley85JOTA}: `Second, it is not clear how minimax and maximin points are related to equilibrium points, except that a joint minimax and maximin point is obviously an equilibrium point.' Moreover, the interpretation as a worst case insurance failed as remarked by Nieuwenhuis \cite[p. 473]{Nieuwenhuis83JOTA}: `Whereas the notions of minimax point and maximin point do not seem to have an easy game-theoretic interpretation, the notion of saddle point has.'

An attempt to resolve this issue by ``brutal force" was given by Ghose and Prasad \cite{GhosePrasad89JOTA}, called Pareto optimal security strategies (POSS). It relies on the idea of providing insurance against the component-wise worst cases, see Ghose \cite{Ghose91JOTA}, Fernandez, Puerto \cite{FernandezPuerto96JOTA}, Fernandez et al \cite{FernandezMonroyPuerto98JOTA} for more details, and leads to very conservative strategies which make little sense in many cases (see Section \ref{sec:rel} below).

The third major difficulty is more mathematical in nature. It boils down to the fact that the infimum and the supremum with respect to vector orders (or even more general preferences) are of no use in most cases. Either they do not exist since the (vector) order is not a lattice order, or even if they do, they produce ``ideal points" which in many cases are not attainable payoffs. The two difficulties described above of generalizing the von Neumann minimax theory to non-complete orders can be attributed to this feature.

Our approach via set relations opens a way out of this dilemma and is strictly motivated by the interpretation of the payoff as multi-dimensional loss for one and gain for the other player. We introduce a new solution concept for  two-player, zero-sum matrix game where main tools are taken from the recently developed theory of set optimization \cite{HamelEtAl15Incoll}, but different ``set orders" as, for instance, in \cite{Maeda15Incoll} are used. Examples illustrate relationships to previous solution concepts and provide evidence that optimal strategies as introduced below indeed should be played if the players wish to achieve ``best" protection against losses, and it is easily seen that they enjoy the interchangeability property. By combining the new concept with Shapley's equilibrium notions, new equilibrium concepts are obtained for which existence theorems are given. Finally, an algorithm for computing such strategies is presented and discussed.

The paper is organized as follows. In the next section, the new solution concept is introduced, and existence theorems for corresponding equilibrium notions are provided. In Section \ref{sec:rel} the new concept is compared with others from the literature. Section \ref{sec:alg} is devoted to an algorithm for computing optimal strategies, and in Section \ref{sec:num} some numerical results are reported. The appendix contains a few crucial concepts from set optimization and also introduces some notation which is used in the paper. We recommend to browse it first in particular to readers who are not familiar with set relations.

\section{A new solution concept} \label{sec:sol}

Let $K \geq 1$ be an integer and $z \in \R^K$ be an element of the $K$-dimensional vector space of column vectors with real components. We write $z = (z_1, z_2, \ldots, z_K)^T$ where the upper $T$ indicates the transpose of the row vector. Let $G = \of{g_{ij}}_{m \times n}$ be an $m \times n$ matrix whose entries 
\[
g_{ij}  = 
	\of{	g^1_{ij}, g^2_{ij}, \ldots, g^K_{ij}}^T \in \R^K
\]
are $K$-dimensional column vectors of real numbers. We interpret $G$ as a loss matrix for the row-choosing player I. Independently, player I and the column-choosing player II select a row $i \in \cb{1, \ldots, m}$ and a column $j \in \cb{1, \ldots, n}$, respectively, which results in player I delivering $g_{ij}$ to player II with the usual convention that negative delivery is something that I receives from II. A possible interpretation of the payoff vectors $g_{ij}$ is, of course, that the $K$ components of $g_{ij}$ denote units of $K$ different assets, so I hands II a portfolio instead of an amount of one particular currency. 

As usual in the theory of finite matrix games, we consider mixed strategies. Player I chooses row $i$ with probability $p_i$ and player II column $j$ with probability $q_j$. The two sets
\[
P = \cb{p \in \R^m_+ \mid \sum_{i=1}^m p_i = 1} \quad \text{and} \quad Q = \cb{q \in \R^n_+ \mid \sum_{j=1}^n q_i = 1}
\]
model the admissible (mixed) strategies of the two players. If player I chooses strategy $p \in P$ and player II $q \in Q$, then the expected (vector) loss of player I is 
\[
v(p, q) = \sum_{i=1}^m\sum_{j=1}^n p_i g_{ij}q_j \in \R^K. 
\]
The set 
\[
V = \cb{v(p, q) \mid p \in P, \; q \in Q}
\]
is in fact a (possibly non-convex) subset of the convex hull of the entries of the matrix $G$.

Different from the scalar case, it is not a priori clear what kind of order should be used for comparing the (expected) payoffs. It would be reasonable to assume that each player has a preference for the payoffs and that these preferences are different from each other--as in Bade \cite{Bade05ET}, but in contrast to references like Nieuwenhuis \cite{Nieuwenhuis83JOTA}, Tanaka \cite{Tanaka94JOTA} where a general partial order in $\R^K$ is used which is the same for both players. The point of view in this note is that it is only known to the players that player I prefers``less loss" and player II prefers``more gain." 

Therefore, it is assumed that both players' decisions are consistent with the partial order $\leq_{\R^K_+}$ generated by the closed convex cone $\R^K_+$, i.e. for $y, z \in \R^K$ the symbol $y  \leq_{\R^K_+} z$ means $z - y \in \R^K_+$, and that no other information about the preferences of the players is available. This is ``Shapley's assumption" \cite[p. 59]{Shapley59NRLQ}: `It is also assumed that the first player wants to increase the components of the vector, and the second player wants to decrease them. Finally, it is assumed that neither player has an a priori opinion concerning the relative importance to himself of the different components.'

Let player I choose a strategy $p\in P$. If for another strategy $p' \in P$ it holds
\[
\forall q \in Q \colon v(p', q) \leq_{\R^K_+} v(p, q),
\]
then clearly $p'$ is better than $p$ for player I. In this case, we write $p' \leq_I p$. In case of
\[
\forall q \in Q \colon v(p', q) = v(p, q),
\]
the strategies $p,p', \in P$ are considered to be equivalent and we write $p' =_I p$. It is natural to determine minimal elements with respect to the order $\leq_I$ and the equivalence relation $=_I$. We say $p \in P$ is {\em $\leq_I$-minimal} if
\[
\of{p' \in P,\; p' \leq_I p} \; \implies\; p' =_I p.
\] 
Likewise, we introduce $\leq_{II}$, $=_{II}$ and {\em $\leq_{II}$-maximality} for player II.

It turns out, however, that the orders $\leq_I$ and $\leq_{II}$ are not ``rich enough"--there are just too many $\leq_I$-minimal and $\leq_{II}$-maximal elements since many cannot be compared with each other.

\begin{example} \label{ex:04} 
In the game given by
\[
G = 
	\of{
	\begin{array}{cc}
	\of{
	\begin{array}{c}
	0 \\ 0
	\end{array}}
	\of{
	\begin{array}{c}
	4 \\ 4
	\end{array}}
	\\[.5cm]
	\begin{array}{cc}
	\of{
	\begin{array}{c}
	3 \\ 1
	\end{array}}
	\of{
	\begin{array}{c}
	1 \\ 3
	\end{array}}
	\end{array}
	\end{array} 
	}
\]
any two strategies $p,p' \in P$ are not comparable with each other with respect to the $\leq_I$-order, and hence they are all $\leq_I$-minimal. Indeed, choosing $\bar q=(1,0)^T$ we have 
\[
 v(p,\bar q) = \sum_{i=1}^m p_i g_{i1} = p_2 \cdot \begin{pmatrix} 3 \\ 1 \end{pmatrix} =
(1-p_1) \cdot \begin{pmatrix} 3 \\ 1 \end{pmatrix},
\]
i.e., the larger $p_1 \in [0,1]$ the smaller with respect to $\leq_{\R_+^K}$ is $v(p,\bar q)$. On the other hand, taking $\hat q = (0,1)^T$ we get
\[
v(p,\hat q) = \sum_{i=1}^m p_i g_{i2} = p_1 \cdot \begin{pmatrix} 4 \\ 4\end{pmatrix} +p_2 \cdot \begin{pmatrix} 1 \\ 3\end{pmatrix},
\]
i.e. the larger $p_1 \in [0,1]$ the larger with respect to $\leq_{\R_+^K}$ is $v(p,\hat q)$. This proves that any two strategies $p,p' \in P$ are not comparable and hence all $p \in P$ are $\leq_I$-minimal. 

On the other hand, consider the strategy $\bar p=(0,1)^T$ for player I. Then
\[
v(\bar p,q) = \sum_{j=1}^n q_j g_{2j} = q_1 \cdot \begin{pmatrix} 3 \\ 1 \end{pmatrix} + q_2 \cdot \begin{pmatrix} 1 \\ 3 \end{pmatrix}.
\]
For any two different elements $q,q' \in Q$, we have $v(\bar p,q) \not\leq_{\R_+^K} v(\bar p,q')$. Thus, elements in $q,q' \in Q$ are never comparable with respect to $\leq_{II}$, which shows that all elements $q \in Q$ are $\leq_{II}$-maximal.  

Even though, based on the order relations $\leq_I$ and $\leq_{II}$, any strategy $p\in P$ is optimal for player I and any $q\in Q$ is optimal for player II, a look at the payoff matrix may convince the reader that it might be better for player I to favor the second row and better for II to have a bias towards the second column.
\end{example}

The above example provokes the question which order relation should be used instead. A minimal requirement certainly is that the order for player I maintains the $\leq_I$-order, but strategies should be comparable more often. If player I picks strategy $p \in P$, then her/his expected payoff belongs to the bounded convex polyhedron
\begin{equation}
\label{EqAvValue1}
v_I\of{p}  = \cb{v(p, q) \mid q \in Q} = \co\cb{\sum_{i=1}^m p_i g_{i1}, \, \ldots, \sum_{i=1}^m p_i g_{in}}.
\end{equation}
For given two strategies $p,p' \in P$, the main idea is to compare the sets $v_I(p)$ and $v_I(p')$ by an appropriate set relation. For some introductory remarks on set relations, the reader is referred to the appendix. We define
\[
v_I(p) \leu v_I(p') \quad :\iff \quad v_I\of{p'} \subseteq v_I\of{p} - \R^K_+.
\]
The following statement is now obvious.
\begin{proposition}
\label{PropSetRelationIa}
Let $p, p' \in P$. If $p' \leq_I p$, then $v_I\of{p'} \leu v_I\of{p}$. 
\end{proposition}

Another aspect is uncovered if one asks which expected payoffs player II can generate if I plays $\bar p \in P$. 
The worst case scenario for player I is a maximal point of $v_I(\bar p)$, i.e. $\bar q \in Q$ satisfying
\[
\of{q \in Q, \; v(\bar p, \bar q) \leq_{\R^K_+} v(\bar p, q)} \; \implies \; v(\bar p, \bar q) = v(\bar p, q).
\]
In this case we write $v(\bar p, \bar q) \in \Max v_I(\bar p)$.

If player II knows the strategy $\bar p$ of player I, then (s)he certainly will pick such a maximal (= non-dominated) strategy. Thus, it might not be necessary to compare all values $\cb{v(\bar p, q) \mid q \in Q}$, but only those which produce maximal elements.

The following result establishes a partial counterpart to Proposition \ref{PropSetRelationIa}.

\begin{proposition}
\label{PropSetRelationIb}
If $p,p' \in P$ with $v_I\of{p'} \leu v_I\of{p}$ and $q \in Q$ with $v(p,q) \in \Max v_I(p)$, then either $v(p', q) \leq_{\R^K_+} v(p, q)$ or $v(p', q), v(p, q)$ are not comparable with respect to $\leq_{\R^K_+}$.
\end{proposition}

{\sc Proof.} Assume $v_I\of{p'} \subseteq v_I\of{p} - \R^K_+$
and $v(p, q) \leq_{\R^K_+} v(p', q)$ for some $v(p,q) \in \Max v_I\of{p}$. Then, by the first statement, there is $\bar q \in Q$ such that $v(p', q) \leq_{\R^K_+} v(p, \bar q)$. Hence
\[
v(p,q) \leq_{\R^K_+} v(p', q) \leq_{\R^K_+}  v(p, \bar q).
\]
Since $v(p,q) \in \Max v_I\of{p}$ this implies $v(p,q) = v(p', q) = v(p, \bar q)$.\pend

\medskip Propositions \ref{PropSetRelationIa} and \ref{PropSetRelationIb} may be interpreted as follows: The transition from the order $\leq_I$ in $P$ to the set order $\leu$ for the comparison of the values $v_I\of{p}$ does not loose comparability, but may absorb some incomparable maximal elements. Thus, the relation $\leu$ produces never more, but often fewer minimal elements: in Example \ref{ex:04}, all strategies $p \in P$ are minimal with respect to the order $\leq_I$, but not all are for $\leu$ (likewise for player II) as is explained in Example \ref{ex:04a} below.

The above discussion almost forces the following set optimization approach to zero-sum matrix games with vector payoff.  Consider the set-valued map defined by
\[
V_I\of{p} := v_I\of{p} - \R^K_+. 
\]
Formally, $V_I$ can (and will) be understood as a function mapping $P$ into $\G(\R^K, -\R^K_+) = \cb{A \subseteq \R^K \mid A =  \cl\co(A - \R^K_+)}$ (see appendix for more details). 

The set $V_I\of{p}$ includes all potential losses for the first player which are less than or equal to $v_I\of{p}$ and thus also absorbs losses generated by ``gifting" something to player II. This is just a version of the standard free disposal condition in economics, see \cite[p. 131]{GreenMasCollelWhinston95Book}. A strategy $p'$ clearly is preferable over $p$ for player I if  $V_I\of{p'} \subseteq V_I\of{p}$, i.e. one can reach the losses in $V_I\of{p'}$ from those in $V_I\of{p}$ by gifting--which usually does not happen. This motivates the following definition.

\begin{definition}
\label{DefSolPI} 
A strategy $\bar p\in P$ is said to be {\em minimal} for player I if there is no $p \in P$ with
\[
V_I\of{p} \subseteq V_I\of{\bar p} \quad \text{and}\quad V_I\of{p} \neq V_I\of{\bar p}.
\]
The set of minimal strategies of player I is denoted by \MIN.	
\end{definition}

In the light of the previous discussion, a minimal strategy is a worst case insurance: The maximal loss, expected to be a maximal point in $v_I(p )$, should be as ``small" as possible which is, in the sense of Definition \ref{DefSolPI}, the case if the set $V_I\of{\bar p}$ is as ``small" as possible. Since $\subseteq$ is, in general, a non-total partial order, this means $V_I\of{\bar p}$ is a minimal element with respect to $\subseteq$.

It is important to note that minimal strategies are independent of the choice of the second player--in contrast to all Nash-type equilibrium notions in the literature  (unless the preference is complete).

\begin{remark}
\label{RemScalarMinimal} If $K = 1$, then $v_I(p) = \cb{v(p, q) \mid q \in Q}$ and
\[
 V_I(p) = v_I(p) - \R_+ = \max_{q \in Q}v(p,q) - \R_+.
\]
Thus, a minimal strategy $\bar p \in P$ satisfies
\[
V_I(\bar p) = \min_{p \in P} \max_{q \in Q} v(p,q) - \R_+ = \max_{q \in Q} v(\bar p,q) - \R_+,
\]
i.e. it is a minimax-strategy for player I.
\end{remark}

A parallel discussion can be done for player II. With
\[
v_{II}\of{q} =  \cb{v(p, q) \mid p \in P} = \co\cb{\sum_{j=1}^n q_j g_{1j}, \, \ldots, \sum_{j=1}^n q_j g_{mj}}
\]
it is clear that player II wants to ``maximize" the set-valued function $q \mapsto v_{II}\of{q}$. The extension of the order $\leq_{II}$ on $Q$ is the set relation $\lel$ for comparing the values of $v_{II}$, which are again convex bounded polyhedra. We define
\[
v_{II}(q') \lel v_{II}(q) \quad :\iff \quad v_{II}\of{q} \subseteq v_{II}\of{q'} + \R^K_+.
\]
Of course, there are analogous statements to Propositions \ref{PropSetRelationIa} and \ref{PropSetRelationIb}. 

A set-valued map is defined by 
\[
V_{II}\of{q} := v_{II}\of{q} + \R^K_+.
\]
Formally, it is a function $V_{II} \colon P \to \G(\R^K, \R^K_+) = \cb{A \subseteq \R^K \mid A =  \cl\co(A + \R^K_+)}$ (see appendix).
It provides the counterpart for $V_I$ with a similar interpretation.

\begin{definition}
\label{DefSolPII} 
A strategy $\bar q \in Q$ is said to be {\em maximal} for player II if there is no $q \in Q$ with
\[
V_{II}\of{q} \supseteq V_{II}\of{\bar q} \quad \text{and} \quad V_{II}\of{q} \neq V_{II}\of{\bar q}.
\]
The set of maximal strategies of player II is denoted by \MAX.	
\end{definition}

Parallel to Remark \ref{RemScalarMinimal}, if $K = 1$ a maximal strategy is a maximizer of the function $q \to \min_{p \in P} v(p,q)$, i.e. a maximin strategy. This shows that for $K=1$ the above concepts boil down to the classic von Neumann approach. An additional (duality) argument is needed to show that the two problems are dual to each other and have the same value.

\begin{example}\label{ex:04a}
Consider the game of Example \ref{ex:04} above. An easy calculation shows that
\[
\MIN = \cb{ p \in P \big|\; 0 \leq p_1 \leq \frac{1}{3}} \; \text{ and } \; \MAX = \cb{ q \in Q \big|\; 0 \leq q_1 \leq \frac{1}{2}}.
\]
The strategy $\hat p = (\frac{2}{3}, \frac{1}{3})^T$ is not minimal for player I; the worst case expected payoff is $\Max v_I(\hat p) = (3, \frac{11}{3})^T$. By playing the minimal strategy $\bar p = (\frac{1}{3}, \frac{2}{3})^T$ with $V_I(\bar p) \subset V_I(\hat p)$ player I can reduce her/his worst case expected payoff to $\Max v_I(\bar p) = (2, \frac{10}{3})^T$. On the other hand, there does not exist another strategy $\bar{\bar p} \in P$ satisfying $\Max v_I(\bar{\bar p}) \subseteq v_I(\bar p) - \R^K_+$, i.e. it is not possible for player I to guarantee a worst case expected payoff strictly better than $(2, \frac{10}{3})^T$ without generating other potential expected payoffs which are not comparable to $(2, \frac{10}{3})^T$ (and hence might be chosen by II).
\end{example}

Another link to set relations should be pointed out. It holds
\[
V_I\of{p} = \bigcup_{q \in Q} \sqb{v(p, q) - \R^K_+} = \sup_{q \in Q} \cb{v(p, q) - \R^K_+}
\]
where the supremum is understood in $(\G(\R^K, -\R^K_+), \subseteq)$, see the appendix for the definition (closure and convex hull can be dropped from the supremum formula as $v(p,\cdot)$ is linear and Q is a convex polyhedron). Since according to Definition \ref{DefSolPI}, player I looks for minimizers of the function $V_I$, her/his problem can be understood as a version of worst case analysis: the worst loss (as described by the supremum) should be minimized. 

Completely parallel to the scalar case, one can look for the infimum of $V_I$, i.e.
\begin{equation}\label{eq:infsup}
	\V_I = \inf_{p \in P}\sup_{q \in Q} \cb{v(p, q) - \R^K_+} = \bigcap_{p\in P}\bigcup_{q \in Q} \sqb{v(p, q) - \R^K_+}.
\end{equation}
This makes sense since $\of{\G(\R^K, -\R^K_+), \subseteq}$ is a complete lattice (see Proposition \ref{PropComLatt} of the appendix). The corresponding problem for player II is to look for the supremum of $V_{II}$ in $\of{\G(\R^K, \R^K_+), \supseteq}$, i.e.
\begin{equation}\label{eq:supinf}
\V_{II} = \sup_{q \in Q}\inf_{p \in P}\cb{v(p, q) + \R^K_+} = \bigcap_{q \in Q}\bigcup_{p\in P} \sqb{v(p, q) + \R^K_+}.
\end{equation}
These two problems are posed in two different image spaces and with respect to different order relations. Therefore, they do not produce a common equilibrium value and cannot be dual--in the sense of linear programming duality--at least not in the same way as the corresponding problems in the scalar case.

In general, the outer infimum in \eqref{eq:infsup} and the outer supremum in \eqref{eq:supinf} are not ``attained'' in a single strategy. Therefore, it cannot be expected that there is a single payoff (vector) which can be considered as the value of the game, and hence there is a multitude of optimal strategies for each  player leading to different (non-comparable) payoffs. This is the reason why infimum and supremum in \eqref{eq:infsup} and \eqref{eq:supinf} are replaced by minimality and maximality notions as introduced in Definitions \ref{DefSolPI} and \ref{DefSolPII}. However, one can show that the sets of minimal and maximal strategies are non-empty and form solutions of the set optimization problem \eqref{eq:infsup} and \eqref{eq:supinf}, respectively, in the sense of \cite[Definition 2.7]{HeydeLoehne11Opt} (compare Definition \ref{DefLatticeSolution} in the appendix). The following theorem provides the essence of the argument.

\begin{theorem} \label{th_ex1} 
For each $(p,q) \in P \times Q$ there exists $(\bar p,\bar q) \in \MIN \times \MAX$ with $V_I(\bar p) \subseteq V_I(p)$ and $V_{II}(\bar q) \subseteq V_{II}(q)$.		
\end{theorem}

{\sc Proof.} This follows from Proposition 5.15 in \cite{HeydeLoehne11Opt}, which states that the domination property (i.e. for every $x \in X$ there exists a minimal point $y \in f[X]:=\cb{f(x) \mid x \in X}$ with $y \leq f(x)$) holds for a function $f \colon X \to Z$, where $X$ is a compact topological space and $(Z,\leq)$ a partially ordered set, whenever $f$ is level-closed. The latter means that for all $z \in Z$ the level sets $\mathcal{L}_f(z) := \cb{x \in X \mid f(x) \leq z}$ are closed. The proof for level-closedness is subject to the following lemma.
\pend

\begin{lemma}
\label{lem_levelclosed} The two functions $V_I \colon P \to \of{\G(\R^K, -\R_+^K), \subseteq}$ and $V_{II} \colon Q \to \of{\G(\R^K, \R_+^K), \supseteq}$ are level-closed.
\end{lemma}

{\sc Proof.} The proof is given for $V_I$ and runs in a similar way for $V_{II}$.

Take $A \in \G(\R^K, -\R^K_+)$ and $\cb{p^\ell}_{\ell = 1, 2, \ldots} \subseteq P$ with $\lim\limits_{\ell \to \infty} p^\ell = \bar p$ in $P$ such that
$V_I(p^\ell) \subseteq  A$ for all $\ell = 1, 2, \ldots$ Then,
\begin{equation}\label{eq_conv_rep}
\forall p \in P \colon	V_I(p)=\co \cb{\sum_{i=1}^m p_i g_{i1},\dots,\sum_{i=1}^m p_i g_{in}} - \R_+^K.
\end{equation}
Thus, for all $\ell = 1, 2, \ldots$,
\[
\forall j=1,\dots, n \colon \sum_{i=1}^m p^\ell_i g_{ij} \in A.
\]
Since $A$ is closed, the same holds for the limit:
\[
\forall j=1,\dots, n \colon \sum_{i=1}^m \bar p_i g_{ij} \in A.
\]
Now, $V_I(\bar p) \subseteq A$ follows from  \eqref{eq_conv_rep} and $A \in \G(\R^K, -\R^K_+)$.  
\pend

\begin{remark}
In the sense of Definition \ref{DefLatticeSolution} below, the set $\MIN$ even is a full solution of problem \eqref{eq:infsup}, and the set $\MAX$ is a full solution of problem \eqref{eq:supinf}.
\end{remark}

Assume that the players pick $(\bar p, \bar q) \in P \times Q$. What could be an incentive for player I to switch to another strategy in $P$? 

First, there is a strategy $\tilde p \in P$ satisfying $\Max v_I(\tilde p) - \R^K_+ \subset \Max v_I(\bar p) - \R^K_+$ (strict inclusion). This means $v_I(\tilde p)- \R^K_+ \subset v_I(\bar p) - \R^K_+$ since, due to the upper domination property (see appendix), $v_I(p) - \R^K_+ = \Max v_I(p) - \R^K_+$ for all $p \in P$. The switch to $\tilde p$ would avoid some potential losses for player I and thus improve her/his worst case estimate (compare Example \ref{ex:04a}).

Secondly, such an incentive is $v(\bar p, \bar q) \not\in \Min v_{II}(\bar q)$, i.e. there is a strategy $p \in P$ such that $v(p, \bar q) \in v_{II}(\bar q)$ with $v( p, \bar q)  \leq_{\R^K_+} v(p, \bar q)$ and $v(p, \bar q)  \neq v(\bar p, \bar q)$. In this case, it makes sense to look for $\hat p \in P$ such that $v(\hat p, \bar q) \in \Min v_{II}(\bar q)$. Such a $\hat p$ always exists since the set $v_{II}(q)$ satisfies the lower domination property (see the appendix for a definition) for all $q \in Q$. 

The first case means that $\bar p$ is not minimal for player I. In the second case $(\bar p, \bar q)$ is not a Shapley equilibrium strategy (compare the following definition) which may happen even if $\bar p$ is minimal (see Example \ref{ex:06}). On the other hand, a strategy can produce a Shapley equilibrium, but not be minimal (see Example \ref{ex:05}). The following two definitions are motivated by these considerations.

\begin{definition}
\label{DefShapleyEq}
A pair $(\bar p, \bar q) \in P \times Q$ is called a Shapley equilibrium if 
\[
v(\bar p, \bar q) \in \Max v_I(\bar p) \cap \Min v_{II}(\bar q).
\]
It is called a strong Shapley equilibrium if 
\[
V_I(\bar p) \cap V_{II}(\bar q) \subseteq \Max v_I(\bar p) \cap \Min v_{II}(\bar q).
\]
\end{definition}

Since
\[
\Max v_I(p) \cap \Min v_{II}(q) \subseteq V_I(p) \cap V_{II}(q)
\]
is always true, the condition for a strong Shapley equilibrium in $(\bar p, \bar q)$ actually means
\[
V_I(\bar p) \cap V_{II}(\bar q) = \Max v_I(\bar p) \cap \Min v_{II}(\bar q).
\]

While Shapley equilibria have been defined in Shapley \cite{Shapley59NRLQ}, strong Shapley equilibria  seem to be a new concept. Clearly, a strong Shapley equilibrium also is a Shapley equilibrium. While the former produces a payoff which cannot be improved by either player with respect to the chosen strategies, the latter  produces payoffs which cannot be improved with respect to the worst case estimate.

The reader is referred to Example \ref{ExCorley3.2} (v) which shows that there are Shapley equilibria which are not strong. Clearly, the feature of being ``strong" can be considered as a refinement of a Shapley equilibrium. However, no additional exogenous or `endogenous parameters' (\cite[page 171]{DeMarcoMorgan07IGTR}) are introduced to the problem which is in contrast to previous approaches.

\begin{remark}
\label{RemShapleyScalarized}
Shapley equilibria can be found by solving non-zero-sum scalar games. More precisely, a pair $(\bar p, \bar q) \in P \times Q$ is a Shapley equilibrium if, and only if, there are $\alpha, \beta \in \Int \R^K_+$ such that $\bar p$ is optimal for player I for the scalar game with the payoff matrix containing the entries $\alpha^Tg_{ij}$, and $\bar q$ is optimal for player II for the game with the matrix $\beta^Tg_{ij}$. This was established in \cite[Theorem]{Shapley59NRLQ}. Note that Shapley denoted such an equilibrium as a Strong Equilibrium Point (SEP).
\end{remark}

As already remarked, minimal/maximal strategies provide worst case payoff estimates for each player independent of the choice of the other. Therefore, minimal/maximal strategies are (trivially) interchangeable. On the other hand, Shapley equilibrium strategies involve both players and finding them is a recursive procedure; interchangeability is violated as already shown by Corley \cite[Example 3.2]{Corley85JOTA}.

\begin{definition}
\label{DefEquilibrium} A pair $(\bar p, \bar q) \in P \times Q$ is called a {\em set relation equilibrium} if $\bar p$ is minimal and $\bar q$ is maximal. 

A set relation equilibrium $(\bar p, \bar q)$ is called a {\em set Shapley equilibrium} if it is also is a Shapley equilibrium.

A set relation equilibrium $(\bar p, \bar q)$ is called a {\em strong set Shapley equilibrium} if it is also is a strong Shapley equilibrium.
\end{definition}

Definition \ref{DefEquilibrium} can be understood as an equilibrium version of the definition of solutions for set optimization problems due to Heyde and L\"ohne in \cite{HeydeLoehne11Opt} (see Definition \ref{DefLatticeSolution} in the appendix). The two features ``being minimal/maximal" and ``attaining an equilibrium value" are no longer equivalent as in the scalar case. 

Example \ref{ExCorley3.2} (v) and Example \ref{ex:01} below show that a set Shapley equilibrium does not need to be strong.

The following theorem ensures the existence of (strong) set Shapley equilibria.

\begin{theorem}
For every zero-sum matrix game with vector payoffs there exists a strong set Shapley equilibrium.	
\end{theorem}

{\sc Proof.} Consider the scalar zero-sum matrix game given by the matrix
\[
\of{\sum_{k=1}^K g^k_{ij}}_{m \times n} = (e^T g_{ij})_{m \times n},
\]
where $e=(1,\dots,1)^T \in \R^K$ and let $(\hat p, \hat q)$ be an equilibrium point for this game, the existence of which follows from linear programming duality. The expected payoff $\nu$ for the scalar game is related to the payoff $v$ of the game with vector payoffs by
\[
\nu(p,q) :=  \sum_{i,j} p_i (e^T g_{ij}) q_j = e^T \sum_{i,j} p_i g_{ij} q_j= e^T v(p,q).
\]
Since $(\hat p, \hat q)$ is an equilibrium for the scalar game, we have
\[
t:= \nu(\hat p, \hat q) = \max_{q \in Q} \nu(\hat p, q) = \min_{p \in P} \nu(p, \hat q).
\]
This can be written as
\[
t = \max_{y \in v_I(\hat p)} e^T y =  \min_{y \in v_{II}(\hat q)} e^T y.
\]
Because of $\max_{y \in -\R^K_+} e^T y  = \min_{y \in \R^K_+} e^T y = 0$, this is equivalent to
\[
t = \max_{y \in V_I(\hat p)} e^T y =  \min_{y \in V_{II}(\hat q)} e^T y.
\]
We conclude
\[
V_I(\hat p) \subseteq H_- := \cb{y \in \R^K \mid e^T y \leq t} \quad \text{and} \quad
 V_{II}(\hat q) \subseteq H_+ := \cb{y \in \R^K \mid e^T y \geq t}.
\]
By Theorem \ref{th_ex1} there exists $(\bar p,\bar q) \in \MIN \times \MAX$ with $V_I(\bar p) \subseteq V_I(\hat p)$ and $V_{II}(\bar q) \subseteq V_{II}(\hat q)$.
 Thus 
\[
V_I(\bar p) \subseteq H_- \qquad \text{and} \qquad V_{II}(\bar q) \subseteq H_+ .
\]
One has $v(\bar p,\bar q) \in V_I(\bar p) \cap V_{II}(\bar q)  \subseteq H_- \cap H_+ =: H.$ 
Hence 
\[
t = e^T v(\bar p, \bar q) = \max_{y \in V_I(\bar p)} e^T y =  \min_{y \in V_{II}(\bar q)} e^T y.
\]
The well-known characterization of a vector minimum (and a vector maximum) by a weighted sum scalarization, see e.g. Zeleny \cite{Zeleny74}, Ehrgott \cite{Ehrgott05}, yields
\[
v(\bar p, \bar q) \in \Max v_I(\hat p)  \cap \Min v_{II}(\hat q).
\]
Thus, $(\bar p,\bar q)$ is an equilibrium. \pend

\medskip The concepts are illustrated by means of the following example which is a version of Corley \cite[Example 3.2]{Corley85JOTA} adapted to our setting.

\begin{figure}
	    \includegraphics[width=0.325\textwidth]{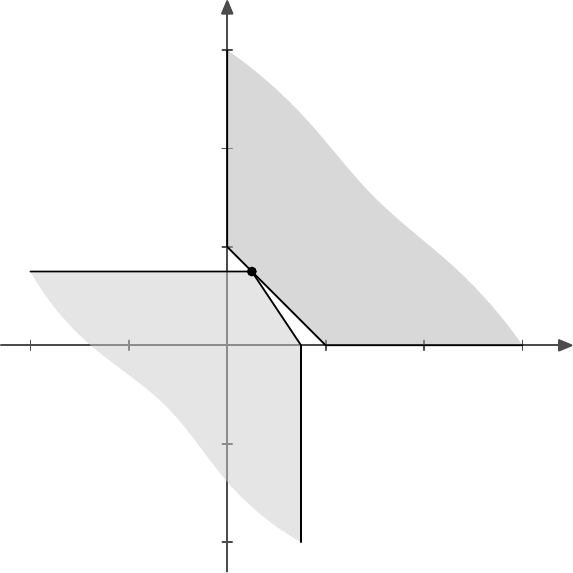}
        \includegraphics[width=0.325\textwidth]{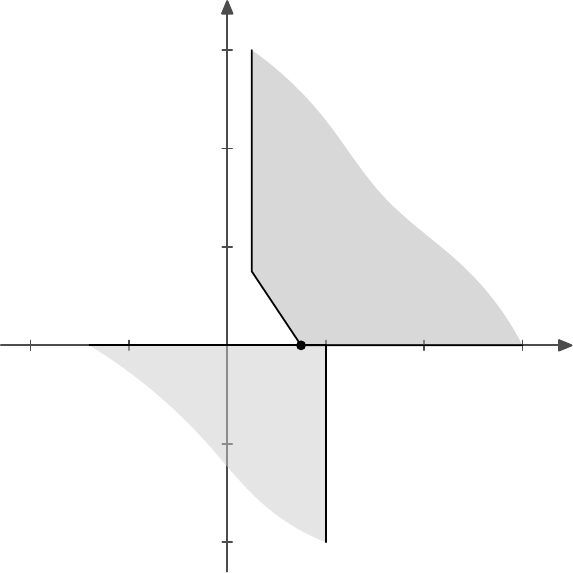}
        \includegraphics[width=0.325\textwidth]{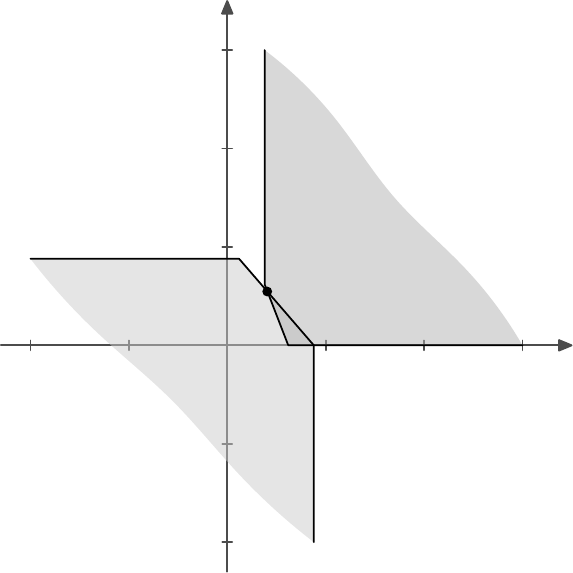}
        
        \caption{The sets $V_I(p)$ and $V_{II}(q)$ for certain strategy pairs $(p,q)$. Left: strong set Shapley equilibrium from Example \ref{ExCorley3.2} (ii); central: the minimal/maximal strategy pair of Example \ref{ExCorley3.2} (iii) is not a Shapley equilibrium; right: the set Shapley equilibrium from Example \ref{ExCorley3.2} (v) is not a strong set Shapley equilibrium.}
        \label{fig:ExCorley3.2}
\end{figure}

\begin{example}
\label{ExCorley3.2} The following facts can be verified for the game given by
\[
G = 
	\of{
	\begin{array}{cc}
	\of{
	\begin{array}{c}
	1 \\ 0
	\end{array}}
	\of{
	\begin{array}{c}
	0 \\ 0
	\end{array}}
	\\[.5cm]
	\begin{array}{cc}
	\of{
	\begin{array}{c}
	0 \\ 1
	\end{array}}
	\of{
	\begin{array}{c}
	1 \\ 0
	\end{array}}
	\end{array}
	\end{array} 
	}.
\]

(i) The set of minimal strategies is $\MIN = \cb{p \in P \mid 0 < p_1 \leq 1}$ and the set of maximal strategies is $\MAX = \cb{q \in Q \mid \frac{1}{2} \leq q_1 \leq 1}$.

(ii) The pairs $(\bar p, \bar q) = \of{(1, 0)^T, (1, 0)^T}$ and $(\hat p, \hat q) = \of{\of{\frac{1}{4}, \frac{3}{4}}^T, \of{\frac{3}{4}, \frac{1}{4}}^T}$ are Shapley equilibria. Moreover, $(\bar p, \bar q)$ is a strong set Shapley equilibrium, see Figure \ref{fig:ExCorley3.2} (left).

(iii) The pair $(\bar p, \hat q)$ is minimal/maximal, but not a Shapley equilibrium, hence not a set Shapley equilibrium. Therefore, there is an incentive for player II to change to $\bar q$ with $v(\bar p, \bar q) = (1, 0)^T$ instead of $v(\bar p, \hat q) = (\frac{3}{4}, 0)^T$, see Figure \ref{fig:ExCorley3.2} (central).

(iv) The pair $(\tilde p, \bar q)$ with $\tilde p = (0,1)^T$ and $v(\tilde p, \bar q) = (0, 1)^T$ is a Shapley equilibrium, but $\tilde p$ is not minimal. There is an incentive for player I to switch. Since player II can generate any payoff in $v_I(\tilde p)$ by an appropriate choice of $q \in Q$, a switch to $p \in P$ with $p_1 \in (0, \frac{1}{2})$ would reduce the degree of freedom for player II considerably, and if $p_1 \in [\frac{1}{2},1]$ is chosen by player I, the choice of player II is even forced (if player II always plays``rational," i.e. maximal points in $v_I(p)$) to generate $v(p, q) = (1- p_1, p_1)^T$ as  this is the only maximal point in $v_I(p)$ for $p_1 \geq \frac{1}{2}$. 

(v) The pair $(p, q)$ with $p = (1/8, 7/8)^T$ and $q = (5/8, 3/8)^T$ is minimal/maximal, a Shapley equilibrium, but not a strong one (see central figure below). Moreover, $p$ is minimal, $\hat q$ is maximal, so $(p, \hat q)$ also is a set Shapley equilibrium, but not a strong set Shapley equilibrium. Figure \ref{fig:ExCorley3.2} (right) shows $V_I(p )$ and $V_{II}(\hat q)$ whose intersection has a non-empty interior and thus includes non-minimal and non-maximal points of $v_I(p )$ and $v_{II}(\hat q)$, respectively.

\end{example}

Again, a strong set Shapley equilibrium can be considered as a refinement of a set Shapley equilibrium.

The question arises if the players have an incentive to switch if they find themselves in a (set) Shapley equilibrium which is not strong. In particular, the players might be tempted to aim for a payoff in $\of{V_I(\bar p) \cap V_{II}(\bar q)}\bs\cb{v(\bar p, \bar q)}$ in case $(\bar p, \bar q)$ is a (set) Shapley equilibrium, but not a strong one. Although it is not possible to achieve a strictly better payoff, a player might be attracted by an alternative payoff which is not comparable to $v(\bar p, \bar q)$ for reasons which are not part of the model so far. The limit is set by the following claim: If $y\in \of{V_I(\bar p) \cap V_{II}(\bar q)}\bs\cb{v(\bar p, \bar q)}$, then $y$ is not comparable to $v(\bar p, \bar q)$. Indeed, if $v(\bar p, \bar q) \leq _{\R^K_+}  y$, then there would be $\bar y \in \Max v_I(\bar p)$ with $v(\bar p, \bar q) \leq _{\R^K_+}  y \leq _{\R^K_+} \bar y$, so  $v(\bar p, \bar q) \not\in \Max v_I(\bar p)$ which contradicts the assumption. A similar argument works for $y \leq _{\R^K_+}  v(\bar p, \bar q)$.

The potential transition from such a payoff $v(\bar p, \bar q)$ to another one in $V_I(\bar p) \cap V_{II}(\bar q)$ might occur due to the `taste' of a player in the sense of Ok et al \cite{OkOrtolevaRiella12Economet} if one interprets the expected gain as her/his utility function.

\section{Relationships to other solution concepts} \label{sec:rel}

\subsection{Minimax and maximin strategies}

First, we will compare minimal and maximal strategies to so-called vector minimax and maximin strategies. Both concepts are in some sense a transfer of worst case estimates from the one-dimensional to the multi-dimensional payoff case.

A strategy $\bar p \in P$ is called a minimax strategy if there is $\bar q \in Q$ such that
\[
v(\bar p, \bar q) \in \Min \bigcup_{p \in P} \WMAX v_I(p).
\]
Vice versa, a strategy $\bar q \in Q$ is called a maximin strategy if there is $\bar p \in P$ such that 
\[
v(\bar p, \bar q) \in \Max \bigcup_{q \in Q} \WMIN v_{II}(q)
\]
Here $\WMIN$ and $\WMAX$ refer to maximal and minimal points with respect to the cone $\{0\} \cup \Int \R^K_+$ which usually are called weakly maximal and weakly minimal, respectively. Minimax and maximin strategies can be understood as ``vector criterion solutions" of the set-valued optimization problems
\[
\text{minimize} \quad p \mapsto \WMAX_w v_I(p) \qquad \text{and} \qquad \text{maximize} \quad q \mapsto \WMIN_w v_{II}(q)
\]
over $p \in P$ and $q \in Q$, respectively. This means, one looks for minimal and maximal points of the union of all function values, compare Hamel et al \cite[p. 80, (III)]{HamelEtAl15Incoll} and the references therein for a discussion of this concept. The above definition can be found in Tanaka \cite{Tanaka94JOTA}, for example. Earlier definitions differ insofar as sometimes minimax strategies include ``outer" weakly minimal and weakly maximal points as in Nieuwenhuis \cite{Nieuwenhuis83JOTA} or drop the ``weak" concept altogether as in Corley \cite{Corley85JOTA}. The latter case is denoted as strong minimax and maximin strategies.

The sets of minimax and maximin payoffs do not coincide in general, and they are also different from the vector saddle point payoffs. This has been observed by Corley \cite{Corley85JOTA} and Nieuwenhuis \cite{Nieuwenhuis83JOTA} (see Section \ref{SecIntro}). Tanaka \cite{Tanaka94JOTA} provides sufficient conditions for the existence of two--in general different--minimax and maximin payoffs.

The following examples show that minimax/maximin strategies on the one hand and minimal/maximal strategies on the other hand form independent concepts in general. 

\begin{example}\label{ex:04b} 
Consider again the game of Example \ref{ex:04} above. 

The set of minimax strategies is $\cb{(0,1)^T, (\frac{1}{3}, \frac{2}{3})^T} \subset \MIN$, the set of maximin strategies is $\cb{q \in Q \mid 0 \leq q_1 \leq \frac{1}{2}} = \MAX$ (this is also the set of all strong maximin strategies).

The set of strong minimax strategies is $P \not\subseteq \MIN$ since $(2,2)^T \in v_I(p)$ for all $p \in P$, and this point belongs to $\Max v_I(p)$ for all $p \in P$ with $0 \leq p_1 < \frac{1}{3}$.
\end{example}

\begin{example}
\label{ExCorley3.2Cont} Consider again the game from Example \ref{ExCorley3.2}. 

The set of all minimax strategies is $\cb{(1, 0)^T}$ which is included in the set of all minimal solutions. The set $\cb{p \in P \mid 0 \leq p_1 < \frac{1}{2}}$ is strongly minimax, and this set is neither included, nor does it include the set of minimal solutions. This also shows the possible ``jump" behavior of the change from minimax to strongly  minimax.

The set of all maximin as well as the set of strongly maximin strategies is $Q$, so $\MAX$ is a strict subset of the set of maximin strategies.
\end{example}

By the way of conclusion, the minimax and maximin concepts from vector optimization do not provide a coherent worst case analysis, i.e. they do not yield bounds for worst possible payoffs. Moreover, they can differ significantly from minimal and maximal solutions. On the other hand, minimal and maximal strategies in the sense of Definition \ref{DefSolPI}, \ref{DefSolPII} have a clear game-theoretic interpretation since they provide worst case bounds for the players.

\subsection{Shapley equilibria}

For $\alpha, \beta \in \R^K_+$ we define
 \[
s_I^\alpha := \inf_{p\in P} \sup_{q\in Q} \sum_{i,\, j} p_i \alpha^T g_{ij} q_j
 	\qquad \text{and} \qquad 
	s_{II}^\beta := \sup_{q\in Q} \inf_{p\in P} \sum_{i,\, j} p_i \beta^T g_{ij} q_j.
 \]
A strategy $\bar p \in P$ is called {\em Shapley solution} for player I if there is $\alpha \in \Int \R^K_+$ such that
 \[
s_I^\alpha = \sup_{q\in Q} \sum_{i,\, j} \bar p_i \alpha^T g_{ij} q_j,
 \]
and $\bar q \in Q$ is called Shapley solution for player II if there is $\beta \in \Int \R^K_+$ such that
\[
s_{II}^\beta = \inf_{p\in P} \sum_{i,\, j} p_i \alpha^T g_{ij} \bar q_j.
\]
A pair $(\bar p, \bar q) \in P \times Q$ of Shapley solutions forms a Shapley equilibrium (compare Remark \ref{RemShapleyScalarized}).

For $A\subseteq \R^K$, we consider the following scalarization functionals which are both versions of the support function of $A$:
 \[
 \sigma^+_A(\alpha) := \sup_{x \in A} \alpha^T x
 	\qquad \text{and} \qquad 
	\sigma^-_A(\beta):=\inf_{x\in A} \beta^T x.
 \]
For a collection $\cb{A_i}_{i \in I}$ of sets in $\R^K$ we have (compare e.g. \cite[Lemma 4.14]{HamelEtAl15Incoll})
\[
\sigma^+_{\bigcup_{i \in I} A_i} = \sup_{i \in I} \sigma^+_{A_i}
	\qquad \text{and} \qquad
	\sigma^-_{\bigcup_{i \in I} A_i} = \inf_{i \in I} \sigma^-_{A_i}.
\]
It follows that
\[
\sigma^+_{V_I(p)}(\alpha) = \sup_{q\in Q} \sum_{i,\, j} p_i \alpha^T g_{ij} q_j
	\qquad \text{and} \qquad
	\sigma^-_{V_{II}(q)}(\beta) = \inf_{p\in P} \sum_{i,\, j} p_i \beta^T g_{ij} q_j.
\]
Hence
\[
s_I^\alpha = \inf_{p\in P} \sigma^+_{V_I(p)}(\alpha)
	\qquad \text{and} \qquad
	s_{II}^\beta = \sup_{q\in Q} \sigma^-_{V_{II}(q)}(\beta).
\]
On the one hand, this shows that Shapley equilibria are related to a scalarization of $V_I(p)$ and $V_{II}(q)$, respectively. On the other hand our optimality concepts can also be characterized in terms of $\sigma^+, \sigma^-$.

\begin{proposition}
\label{PropMinMaxScalar}
A strategy $\bar p \in P$ is minimal for player I if, and only if,
\[
\of{p \in P, \; V_I(p) \neq V_I(\bar p)} \; \Rightarrow  \of{\exists \alpha \in \R^K_+ \colon \sigma^+_{V_I(p)}(\alpha) > \sigma^+_{V_I(\bar p)}(\alpha)}.
\]
Likewise, a strategy $\bar q \in Q$ is maximal for player II if and only if
\[
\of{q \in Q, \; V_{II}(q) \neq V_{II}(\bar q)} \; \Rightarrow \; \of{\exists \beta \in \R^K_+ \colon \sigma^-_{V_{II}(q)}(\beta) < \sigma^-_{V_{II}(\bar q)}(\beta)}.
\]
\end{proposition}
{\sc Proof.} This follows from the equivalence
	$$ A \curlyeqprec B \; \iff \;\sigma^+_A \leq \sigma^+_B,$$
	for $\sigma^+_A$ defined on $\R^K_+$.
	Likewise, we have 
	$$ A \preccurlyeq B \;\iff \;\sigma^-_A \leq \sigma^-_B$$
	for $\sigma^-_A$ defined on $\R^K_+$.\pend

The next two examples show that our optimality notions are independent from the Shapley optimality concept.

\begin{example} \label{ex:05}
If $\bar p \in P$ is a Shapley solution for player I, then it is not necessarily minimal. Indeed, if 
\[
G =  \begin{pmatrix}
		\begin{pmatrix}
			0 \\ 0
		\end{pmatrix} &
		\begin{pmatrix}
			0 \\ 0
		\end{pmatrix}\\[.5cm]
		\begin{pmatrix}
			1 \\ -1
		\end{pmatrix}&
		\begin{pmatrix}
			-1 \\ 1
		\end{pmatrix}
	\end{pmatrix},
\]
then $\bar p = (0,1)^T$ is a Shapley solution for player I (with respect to $\alpha = (\frac{1}{2},\frac{1}{2})$), but it is not minimal in the sense of Definition \ref{DefSolPI}, since $V_I(1,0) \subsetneq V_I(0,1)$. In fact,
\[
V_I(1,0) = \cb{\begin{pmatrix} 0 \\0 \end{pmatrix}} - \R^2_+
	\qquad \text{and} \qquad
	V_I(0,1) = \co\cb{\begin{pmatrix} 1 \\-1 \end{pmatrix},\begin{pmatrix} -1 \\1 \end{pmatrix}} - \R^2_+.
\]
\end{example}

\begin{example} \label{ex:06}
If $\bar p \in P$ is minimal for player I in the sense of Definition \ref{DefSolPI}, it is not necessarily a Shapley solution.
Indeed, if
\[
G =  \begin{pmatrix}
		\begin{pmatrix}
			0 \\ 0
		\end{pmatrix} &
		\begin{pmatrix}
			3 \\ -3
		\end{pmatrix}\\[.5cm]
		\begin{pmatrix}
			-3 \\ 3
		\end{pmatrix}&
		\begin{pmatrix}
			0 \\ 0
		\end{pmatrix}\\[.5cm]
		\begin{pmatrix}
			1 \\ 1
		\end{pmatrix}&
		\begin{pmatrix}
			1 \\ 1
		\end{pmatrix}
	\end{pmatrix},
\]
then the strategy $\bar p = (0,0,1)^T$ is minimal, but not Shapley: Consider the sets
\[
V_I(1,0,0) = \co\cb{\begin{pmatrix} 0 \\0 \end{pmatrix},\begin{pmatrix} 3 \\-3 \end{pmatrix}}-\R^2_+,\quad V_I(0,1,0) = \co\cb{\begin{pmatrix} -3 \\ 3 \end{pmatrix},\begin{pmatrix} 0 \\ 0 \end{pmatrix}}-\R^2_+,
\]
as well as
\[
V_I(0,0,1) = \cb{\begin{pmatrix} 1 \\ 1 \end{pmatrix}}-\R^2_+
\]
and take into account the above scalarization results. It is easy to see that $\bar p$ is not a Shapley solution. Assume that $\bar p$ is not minimal. Then there is some $p \in P$ with $V_I(p) \subseteq V_I(\bar p)$ and $V_I(p) \neq V_I(\bar p)$. Then we have
\[
\forall q \in Q:\quad p_1 \begin{pmatrix} 3\\-3 \end{pmatrix} q_2 + p_2 \begin{pmatrix} -3\\ 3\end{pmatrix} q_1 + p_3 \begin{pmatrix} 1\\ 1\end{pmatrix} \in \cb{\begin{pmatrix} 1\\ 1\end{pmatrix} }- \R^2_+.
\]
Since $q_1 + q_2 = 1$ we obtain
\begin{align*}
(-3 p_1 - 3 p_2) q_1 & \leq 1 - p_3 -3 p_1, \\
(3 p_1 + 3 p_2) q_1 & \leq 1 - p_3  + 3 p_1
\end{align*}
for all $q_1 \in [0,1]$. This system is satisfied if and only if the first inequality holds for $q_1=0$ and the second one for $q_1 = 1$. Thus, it is equivalent to
\begin{align*} 
p_3 + 3 p_1 & \leq 1, \\
p_3 + 3 p_2 & \leq 1.
\end{align*}
It follows that $p = (0,0,1)^T$, which shows that $\bar p$ is minimal.
\end{example}

\subsection{Maeda's bi-matrix games with set payoffs}

Maeda \cite{Maeda15Incoll} introduced different concepts of ``equilibrium points'' for bi-matrix games with set payoffs. If one specializes these concepts to  the setting of this note, one may see that Maeda's {\em maximal Nash equilibrium} is closely related to, but different from the optimality concepts introduced and motivated in the previous section. We shortly recall some definitions and results from \cite{Maeda15Incoll} using our notation and our setting of assumptions.

A bi-matrix game with set payoffs is defined by an $m \times n$ matrix the entries of which are pairs $(A_{ij},B_{ij})$ of sets $A_{ij},B_{ij} \subseteq \R^K$. Mixed strategies are considered which results for player I in an expected payoff
$$ v_I(p,q) = \sum_{i=1}^m\sum_{j=1}^n p_i A_{ij} q_{j}  \subseteq \R^K$$
and for player II in
$$ v_{II}(p,q) = \sum_{i=1}^m\sum_{j=1}^n p_i B_{ij} q_{j} \subseteq \R^K.$$
Both players maximize their expected payoffs. By considering the singleton sets
$$ A_{ij}= \cb{-g_{ij}} \qquad B_{ij} = \cb{g_{ij}}$$
we obtain a zero-sum game with vector payoff, as considered in this article.

Maeda \cite{Maeda15Incoll} argues that each player can choose a set relation. In particular, the set relation $\lel$ (``L-type'') and $\leu$ (``U-type'') as introduced in Section \ref{sec:sol} are suggested as candidates for preferences of the players. Moreover, the ``LU-type" relation, which is defined by the requirement that both $\lel$ and $\leu$ are satisfied, is suggested. The definitions of various types of equilibrium points as well as corresponding existence results are given in \cite[p.  320]{Maeda15Incoll} under the assumption that `both players $[\ldots]$ are LU type and this is a common knowledge for the players.'

In \cite[Definition 4.1]{Maeda15Incoll}, a {\em Nash equilibrium strategy} is introduced. In our setting this is a pair $(p^*,q^*) \in P\times Q$ such that
\[
\forall p \in P: p^* \leq_I p \quad \text{and} \quad \forall q \in Q: q \leq_{II} q^*.
\]
The existence of a Nash equilibrium strategy was shown under a very strong assumption, which is not fulfilled, for instance, in Example \ref{ex:04}. We have shown that any two strategies for I are not comparable (likewise for II) and hence there is no Nash equilibrium strategy. This problem has been also addressed in \cite[Example 4.2]{Maeda15Incoll}, and for this reason other concepts have been introduced. A {\em maximal Nash equilibrium strategy} (compare \cite[Definition 4.2]{Maeda15Incoll}) is (in our setting) a pair $(p^*,q^*)$ such that $v_{I}(p^*)$ is minimal and $v_{II}(q^*)$ is maximal with respect to a set relation. But as mentioned above, this set relation is supposed to be of LU-type, whereas in our concept player I has to use $\leu$ (U-type) and player II has to use $\lel$ (L-type). This choice is basically forced by the interpretation of the payoff as loss for player I and gain for player II, and it is also of practical relevance as the following example shows.

\begin{example}
\label{ExZeroRow}
	Consider the game
	\[
	G =  \begin{pmatrix}
			\begin{pmatrix}
				2\\ -1
			\end{pmatrix} &
			\begin{pmatrix}
				-1 \\ 2
			\end{pmatrix}\\[.5cm]			
			\begin{pmatrix}
				0 \\ 0
			\end{pmatrix}&
			\begin{pmatrix}
				0 \\ 0
			\end{pmatrix}	
		\end{pmatrix}.
	\]
	For $p \in P$ and $q \in Q$ we have
	\[ 
	v_{I}(p) = \co\cb{ p_1 \begin{pmatrix}2 \\ -1\end{pmatrix}, p_1 \begin{pmatrix} -1\\2\end{pmatrix} },
	\]	
	\[ 
	v_{II}(q) = \co\cb{ \begin{pmatrix} 0 \\ 0\end{pmatrix}, q_1 \begin{pmatrix} 2\\ -1\end{pmatrix} + q_2 \begin{pmatrix} -1\\ 2\end{pmatrix}}
	\]
	One can easily verify that for any two different $p,p' \in P$, the sets $v_I(p), v_I(p')$ are not comparable with respect to $\lel$. But using the relation $\leu$, for $p^*=(0,1)^T$ we have
	\[
	\forall p \in P:\quad v_I(p^*) \leu v_I(p) \quad \text{and} \quad \forall p \in P\setminus\cb{p^*}:\quad v_I(p^*) \leu v_I(p) 
	\]
	On the other hand, for any two different $q,q' \in Q$, the sets $v_{II}(q), v_{II}(q')$ are not comparable with respect to $\leu$. For the relation $\lel$, setting
	$Q^*= \cb{q \in Q \mid \frac{1}{3} \leq q_1 \leq \frac{2}{3}}$ we have
	\[
	\forall q,q' \in Q^*:\quad v_{II}(q) \lel v_{II}(q') \quad \text{and} \quad v_{II}(q') \lel v_{II}(q).
	\]
	Moreover we have
	\[
	\forall q \in Q,\; \forall q^* \in Q^*:\quad v_{II}(q) \lel v_{II}(q^*).
	\]
We conclude the following for player I: While with our approach only $p^*=(0,1)^T$ is optimal, using the LU-type ordering yields that all $p \in P$ are optimal. For player II it is similar: For our approach the strategies in $Q^*$ are optimal, but using the LU-type ordering all strategies $q \in Q$ are so. 
	
We conclude that using the LU-type ordering may not provide any information to the players. It is, however, quite obvious that a loss averse player I should prefer to choose the second row. Not so obvious, but equally motivated, is that the strategies in $Q^*$ are preferable for player II.
\end{example}

The example shows that existence results for optimality notions of games with vector and set payoff are not sufficient for their justification. It is rather important to ensure that not too many optimal strategies exist and that the concepts come along with a clear motivation and interpretation.

\subsection{Pareto optimal security strategies (POSS)}

Assume that player I picks strategy $p \in P$. The set
\[
W_I(p ) := \bigcap_{q \in Q} \sqb{v(p, q) + \R^K_+}
\]
contains all expected losses which player I can suffer by choosing strategy $p$ independently of player II's choice. The addition of the convex cone $\R^K_+$ reflects the fact that player I always can `gift something" to II. It makes sense to make this set as ``big" as possible, i.e. include as many potential losses as possible since then the chance that there are ``small ones" among them is bigger. Thus, it makes sense to look for 
\[
\bigcup_{p \in P}\bigcap_{q \in Q} \sqb{v(p, q) + \R^K_+}.
\]
Up to a closure and a convex hull, this expression coincides with
\[
\inf_{p \in P}\sup_{q \in Q} \sqb{v(p, q) + \R^K_+}
\]
where $\inf$ and $\sup$ are understood in $\of{\G(\R^K, \R^K_+), \supseteq}$ (see appendix). The expression $\bigcap_{q \in Q} \sqb{v(p,q) + \R_+^K}$ is not changed if $Q$ is replaced by the set if its vertices, hence we get
\begin{equation*}\label{eq:poss1}
\bigcup_{p \in P}\bigcap_{q \in Q} \sqb{v(p, q) + \R^K_+} = \bigcup_{p \in P}\bigcap_{j \in \cb{1, \ldots, n}} \sqb{\sum_{i=1}^m p_i g_{ij} + \R^K_+}.
\end{equation*}
We also conclude that the set $\cb{(z, p) \in \R^K \times P \mid z \in W_I(p)}$ is closed and (polyhedral) convex. Moreover, $P$ is compact. Now it is easy to see that the closed convex hull in the definition of the infimum can be omitted here, i.e., we have
\[
\inf_{p \in P}\sup_{q \in Q} \sqb{v(p, q) + \R^K_+}= \bigcup_{p \in P}\bigcap_{q \in Q} \sqb{v(p, q) + \R^K_+}.
\]

If, as usual, player I tries to minimize her/his maximal expected loss (s)he is led to the following $\G(\R^K, \R^K_+)$-valued problem:
Find
\[
\W_I = \inf_{p \in P} \sup_{j \in \cb{1, \ldots, n}} \sqb{\sum_{i=1}^m p_i g_{ij} + \R^K_+} =
	\bigcup_{p \in P} \bigcap_{j \in \cb{1, \ldots, n}} \sqb{\sum_{i=1}^m p_i g_{ij} + \R^K_+}.
\]

\begin{definition}
\label{DefPOSSI} A strategy $\bar p \in P$ is called a Pareto optimal security strategy (POSS) for player I if there is no $p \in P$ satisfying
\[
W_I\of{p} \supseteq W_I\of{\bar p} \quad \text{and} \quad W_I\of{p} \neq W_I\of{\bar p}.
\]
The set of POSS for player I is denoted by $\text{POSS}(I)$.
\end{definition}

Player II proceeds in a similar way. The set
\[
W_{II}\of{q} := \bigcap_{p \in P} \sqb{v(p, q) - \R^K_+} \in \mathcal G\of{\R^K, -\R^K_+}
\]
includes all potential gains for her/him including those obtained by ``giving up something for free," and this set should be ``as big as possible." So, player II is faced with the problem to find
\[
\W_{II} = \sup_{q \in Q} \inf_{i \in \cb{1, \ldots, m}} \sqb{\sum_{j=1}^n q_j g_{ij} - \R^K_+} =
	\bigcup_{q \in Q} \bigcap_{i \in \cb{1, \ldots, m}} \sqb{\sum_{j=1}^n q_j g_{ij} - \R^K_+}.
\]

\begin{definition}
\label{DefPOSSII} A strategy $\bar q \in Q$ is called a Pareto optimal security strategy (POSS) for player II if there is no $q \in Q$ satisfying
\[
W_{II}\of{q} \subseteq W_{II}\of{\bar q} \quad \text{and} \quad W_{II}\of{q} \neq W_{II}\of{\bar q}.
\]
The set of POSS for player II is denoted by $\text{POSS}(II)$.
\end{definition}

The previous two definitions are versions of Definition 4.1 in Ghose, Prasad \cite{GhosePrasad89JOTA} adopted to our setting. The following results are well-known, see Fernandez, Puerto \cite[Theorem 3.1]{FernandezPuerto96JOTA}.
  
\begin{proposition}
Define the two sets
\begin{align*}
S_I & = \cb{(p, y) \in \R^m \times \R^K \mid  y \geq \displaystyle\sum_{i=1}^m p_ig_{ij},\; j=1,\ldots,n, \; p \geq 0,\; e^T p = 1} \\
S_{II} & = \cb{(p, y) \in \R^m \times \R^K \mid  y \leq \displaystyle\sum_{j=1}^n g_{ij}q_j,\; i=1,\ldots,m, \; q \geq 0,\; e^T q = 1} 
\end{align*}
Then, 
\[
\W_I = \cb{y \mid (p, y) \in S_I} + \R^K_+ \quad \text{and} \quad \W_{II} = \cb{y \mid (p, y) \in S_{II}} - \R^K_+.
\]
\end{proposition} 
 
The result means that Pareto optimal security strategies as well as the sets $\W_I,  \W_{II}$ can be obtained by solving two linear multi-criteria optimization problems (MLOP). This is important for computational approaches. Moreover, the POSS approach is related to the concepts introduced in Section \ref{sec:sol} as follows.

\begin{proposition}
\label{PropWeakDual}
It holds
\begin{align*}
\forall p \in P, \; \forall q \in Q \colon W_{II}\of{q} \subseteq V_I\of{p}, \\
\forall p \in P, \; \forall q \in Q \colon V_{II}\of{p} \supseteq W_I\of{q}.
\end{align*}
Moreover, $\V_{II}  \supseteq \W_I$ and $\W_{II} \subseteq \V_I$.
\end{proposition}
 
{\sc Proof.} Everything is immediate from the definitions. \pend

\medskip  The next result shows that, as a rule, optimal strategies are not worse than POSS with respect to their payoffs, compare also Example \ref{ex:01} below. Therefore, one may guess that optimal strategies even lead to better worst case estimates for the expected payoff, and this is indeed the case for many examples.
 
\begin{theorem} \label{ThmMinMaxPOSS} It holds
\begin{align*}
\forall p \in \MIN & :  V_I(p) \cap \of{\W_I+ \R^K_+\bs\{0\}} = \emptyset, \\
\forall q \in  \MAX & : V_{II}(q) \cap \of{\W_{II} - \R^K_+\bs\{0\}} = \emptyset.
\end{align*}
\end{theorem}

{\sc Proof.} Fix $\bar p \in P$. According to \eqref{EqAvValue1}, one has
\[
V_I(\bar p) = \co\cb{\sum_{i=1}^m \bar p_i g_{i1}, \, \ldots, \sum_{i=1}^m \bar p_i g_{in}} - \R^K_+.
\]
Assume there is 	
\[
z \in V_I(\bar p) \cap \of{\W_I + \R^K_+\bs\{0\}}.
\]
Then there exist $p \in P$ and $c \in \R^K_+\bs\{0\}$ such that both is satisfied $z \in V_I(\bar p)$ and
\[
z - c \in \bigcap_{j \in \cb{1, \ldots, n}} \sqb{\sum_{i=1}^m p_i g_{ij} + \R^K_+}.
\]
This implies
\begin{multline*}
V_I(p) = \co \cb{\sum_{i=1}^m p_i g_{i1}, \, \ldots, \sum_{i=1}^m p_i g_{in}} -\R^K_+ \subseteq \cb{z-c} - \R^K_+ \subsetneq \cb{z} - \R^K_+ \subseteq V_I(\bar p).
\end{multline*}
Thus $\bar p$ is not minimal for player I, which proves the first claim. The second statement can be shown analogously. \pend

By the way of conclusion, POSS can also be obtained by a set optimization approach, but are different from minimal/maximal strategies in general. Example \ref{ex:01} below illustrates the results of this section and shows that Pareto optimal security strategies are very often too conservative as a worst case estimate.

\section{How to compute optimal strategies} \label{sec:alg}
	
The procedure for player I is based on optimality tests for strategies $p \in P$. As there are infinitely many such strategies, one has to choose a finite subset $\bar P$ of $P$ first. The finite family $\T=\cb{V_I(p) \mid p \in \bar P \cap \MIN}$ of payoff sets with respect to optimal strategies is presented to player I who acts as a decision maker. Assuming that the finitely many optimal strategies found in this way provide a good representation of all optimal strategies, player I can use all the information of the finite family $\T$ of polyhedral convex sets to select one optimal strategy $p^*$ to play.

The question to be raised is how such a minimality test can be implemented. Assume we want to test whether or not some $\bar p \in P$ belongs to $\MIN$. We have
\[
V_I(\bar p) = \cb{\sum_{i=1}^m\sum_{j=1}^n \bar p_i g_{ij} q_j \big|\; q \geq 0, \; e^T q = 1}
\]
where $e=(1,\dots,1)^T \in \R^n$. In the first step we compute an {\em H-representation} of the convex polyhedron $V_I(\bar p)$, that is, we compute $H \in \R^{s \times m}$, $h \in \R^s$ such that
\[
V_I(\bar p) = \cb{y \in \R^K \mid H y \geq h}.
\] 
Secondly, we describe the condition $V_I(p) \subseteq V_I(\bar p)$ by linear inequalities with respect to variables $p \in P$.
\begin{proposition}\label{prop:alg1}
	For any $p \in P$, the following is equivalent:
	\begin{enumerate}[(i)]
		\item $V_I(p) \subseteq V_I(\bar p)$
		\item $p$ is feasible for the system of linear inequalities
		\begin{equation}\label{eq:alg1}
			H \cdot \of{\sum_{i=1}^m p_i g_{ij}} \geq h, \quad j=1,\dots,n.			
		\end{equation} 
	\end{enumerate}	
\end{proposition}
{\sc Proof.}  We have
$$ V_I(p) = \co \cb{\sum_{i=1}^m p_i g_{i1}, \dots, \sum_{i=1}^m p_i g_{in}} - \R_+^K$$
for all $p \in P$.
Thus, (i) is satisfied if and only if all points occurring in the convex hull expression belong to $V_I(\bar p)$.\pend

We intend to test whether or not there is $p \in P$ such that $V_I(p) \subseteq V_I(\bar p)$ and $V_I(p) \neq V_I(\bar p)$. This will be done by solving a linear program. We already know that $p \in P$ and  $V_I(p) \subseteq V_I(\bar p)$ can be described by linear inequalities. It remains to find a linear inequality description of the condition $V_I(p) \neq V_I(\bar p)$. To this end we compute the vertices $y^1,\dots, y^r$ of the polyhedral convex set $V_I(\bar p)$.  

\begin{proposition}\label{prop:alg2}
Let $p,\bar p \in P$ and $V_I(p) \subseteq V_I(\bar p)$. Then the following is equivalent:
\begin{enumerate}[(i)]
	\item $V_I(p) \neq V_I(\bar p)$
	\item There exists a vertex $y^i$ of $V_I(\bar p)$ such that $y^i \not\in V_I(p)$.
\end{enumerate}	
\end{proposition}
{\sc Proof.} Obviously, (ii) implies (i). Assume that (ii) is not true, i.e., all vertices of $V_I(\bar p)$ belong to $V_I(p)$. Since $V_I(p) - \R_+^K = V_I(p)$ and 
$$   V_I(\bar p) = \co \cb{y^1,\dots, y^r} - \R_+^K,$$ 
we obtain $V_I(p) \supseteq V_I(\bar p)$ and hence  $V_I(p) = V_I(\bar p)$.\pend

For each vertex $y^\ell$ $(\ell=1,\dots,r)$ of $V_I(\bar p)$ we now compute a supporting hyperplane
\[
H^\ell = \cb{y \in \R^K\mid (c^\ell)^T y = \gamma_\ell} 
\]
to $V_I(\bar p)$ with the property 
\begin{equation}\label{eq:alg3}
V_I(\bar p) \cap H^\ell= \cb{y^\ell}.
\end{equation}
For the representation of $H^\ell$ we assume that
\[
	\forall y \in V_I(\bar p):\; (c^\ell)^T y \geq \gamma_\ell.
\]
By similar arguments as in the proof of Proposition \ref{prop:alg1} we conclude from $V_I(p)\subseteq V_I(\bar p)$ that
\begin{equation}\label{eq:alg0}
	\forall j=1,\dots,n:\; (c^\ell)^T \sum_{i=1}^m p_i g_{ij} \geq \gamma_\ell.
\end{equation}
If for some $p \in P$ with $V_I(p)\subseteq V_I(\bar p)$, some $\ell \in \cb{1,\dots, r}$ and some $\varepsilon_\ell > 0$ the linear system 
\begin{equation}\label{eq:alg2}	
	\forall j=1,\dots,n:\; (c^\ell)^T \sum_{i=1}^m p_i g_{ij} \geq \gamma_\ell + \varepsilon_\ell
\end{equation}
is satisfied, then the vertex $y^\ell$ of  $V_I(\bar p)$ does not belong to $V_I(p)$. Using Proposition \ref{prop:alg2} we conclude $V_I(p)\neq V_I(\bar p)$. This means that $\bar p \not\in \MIN$.

Vice versa, let \eqref{eq:alg2} be violated for all $p \in P$ with $V_I(p)\subseteq V_I(\bar p)$, all $\ell \in \cb{1,\dots, r}$ and all $\varepsilon > 0$.  By \eqref{eq:alg0}, for every $p \in P$ with $V_I(p)\subseteq V_I(\bar p)$ and every vertex $y^\ell$ of  $V_I(\bar p)$ there exists $j_{(P,\ell)} \in \cb{1,\dots,n}$ such that
\[
	(c^\ell)^T \sum_{i=1}^m p_i g_{ij_{(P,\ell)}} = \gamma_\ell.
\]
By \eqref{eq:alg3}, this means that, for some fixed $p$, any vertex $y^\ell$ of $V_I(\bar p)$ coincides with some point of $V_I(p)$. In this situation we have $V_I(p) = V_I(\bar p)$ for all $p \in P$ with $V_I(p)\subseteq V_I(\bar p)$ and hence $\bar p \in \MIN$.

The considerations above can be summarized as follows.

\begin{theorem}
A strategy $\bar p \in P$ for player I is a minimal strategy in the sense of Definition \ref{DefSolPI} if, and only if, the following linear program has the optimal value zero: 
\begin{equation}
	\max_{p,\varepsilon} \sum_{\ell=1}^r\varepsilon_\ell  \; \text{ subject to } \\ \eqref{eq:alg1},\; \eqref{eq:alg2} \text{ and }  p \geq 0, \; e^T p = 1
\end{equation}	
where $e=(1,\dots,1)^T \in \R^m$.	
\end{theorem}

We close this section with some remarks on implementation details.

\begin{remark} \label{rem1} An H-representation as well as the vertices of $V_I(\bar p)$ can be computed by solving a multiple objective linear program (MOLP). The vertices are obtained from the primal problem and an H-representation from the dual problem. Moreover, the hyperplanes $H^\ell$ can be obtained in this way. Only one MOLP needs to be solved to test some $\bar p \in P$ for optimality. For more details, the reader is referred, for instance, to \cite{LoeWei16}.	
\end{remark}

\begin{remark} \label{rem2} Theorem \ref{ThmMinMaxPOSS} can be used to sort out some (but not all) non-optimal strategies. To this end we need an inequality representation of $\W_I$ (POSS payoffs), say
$$ \W_I = \cb{y \in \R^K \mid A y \geq a},$$	
which can be computed by solving a multiple objective linear program. This has to be done only once.
For some $\bar p$, consider the following system of linear inequalities with variables $y\in \R^K$ and $q \in \R^n$:
\begin{equation}\label{eq:rem1}
	q \geq 0,\;e^T q = 1,\;  y \leq \sum_{j=1}^n \sum_{i=1}^m \bar p_i g_{ij} q_j, \;  A (y - \varepsilon e) \geq a.
\end{equation}
The (small) parameter $\varepsilon > 0$ is used to take into account the exclusion of zero in the term $\R^K_+\bs\cb{0}$ in the formula of Theorem \ref{ThmMinMaxPOSS}. If the linear system \eqref{eq:rem1} is feasible, then $\bar p \not\in \MIN$. Feasibility of \eqref{eq:rem1} can be verified by solving an LP. This LP is smaller than the LP in the optimality test. Moreover, it is not necessary to compute an H-representation and the vertices of $V_I(\bar p)$. Therefore it can be more efficient to sort out some non-optimal strategies in a first step and to execute minimality tests only for the remaining strategies.
\end{remark}

Of course, the algorithm for player II is completely analogous. Our procedure to test whether a pair $(p,q) \in \MIN \times \MAX$ is a set Shapley equilibrium uses H-representations of $V_I(p)$ and $V_{II}(q)$, which already have been computed in the optimality tests. To decide whether or not $v(p,q)$ belongs to $\Max v_I(p)$, consider those inequalities of the H-representation of $V_I(p)$ in which equality holds for at the point $v(p,q)$. Let $a$ be the sum of the outer normals of the corresponding supporting hyperplanes. Then $v(p,q) \in \Max v_I(p)$ if, and only if, $a \in \Int \R^K_+$. Likewise, $v(p,q) \in \Min v_{II}(q)$ can be checked.

In order to test whether or not a set Shapley equilibrium $(p,q)$ even is a strong set Shapley equilibrium, we note that the condition 
$$ V_I(p) \cap V_{II}(q) \subseteq \Max V_I(p) \cap \Min V_{II}(q)$$
is equivalent to
\begin{equation} \label{num_strong}
	\forall \varepsilon > 0,\; \forall k \in \cb{1,\dots,K}: \; V_I(p) \cap (V_{II}(q)+\varepsilon e^k) = \emptyset
\end{equation}
where $e^k$ denotes the $k$-th unit vector in $\R^K$. Again we use the H-representations $V_I(p)=\cb{y \in \R^K \mid A y \geq a}$ and $V_{II}(q)=\cb{y \in \R^K \mid B y \geq b}$ and we consider the linear program
\begin{equation} \label{lp_strong}
	 \max e^T t \;\text{ subject to }\; A y \geq a,\; B (y - t) \geq b,\; t \geq 0.
\end{equation}
This LP is bounded (as $V_{II}(q) = v_{II}(q) + \R^K_+$ and $v_{II}(q)$ is bounded) and the point $(y,t)=(v(p,q),0)$ is feasible. The optimal value of \eqref{lp_strong} is zero if, and only if, \eqref{num_strong} holds.

\section{Numerical results} \label{sec:num}

The algorithm of the previous section has been implemented with GNU Octave, version 4.2.
We used the VLP-Solver Bensolve version 2.0.1 \cite{bensolve, LoeWei16} linked against the GLPK library version 4.6 to compute an H-representation and the vertices of $V_I(\bar p)$, compare Remark \ref{rem1}. All computations were run on a computer with Intel\textregistered\ Core\texttrademark\ M CPU with 1.2 GHz and 8GB of RAM . 

To compute a finite representation of the set $\MIN$, the set $P$ is discretized with a stepsize $t>0$, i.e. optimality is tested over the finite set
$$ \bar P = P \cap \cb{y \in \R^K \mid \forall k \in \cb{1,\dots,K},\; \exists z_k \in \Z:\; y_k=z_k t }$$
and likewise for player II.

\begin{example}\label{ex:01} Let us consider the game
\[ 
G =  \begin{pmatrix}
		\begin{pmatrix}
			5 \\ 0
		\end{pmatrix} &
		\begin{pmatrix}
			-1 \\ -5
		\end{pmatrix}&
		\begin{pmatrix}
			4 \\ -4
		\end{pmatrix}\\[.5cm]
		\begin{pmatrix}
			2 \\ -2
		\end{pmatrix}&
		\begin{pmatrix}
			2 \\ -7
		\end{pmatrix}&
		\begin{pmatrix}
			2 \\ 2
		\end{pmatrix}\\[.5cm]
		\begin{pmatrix}
			0 \\ -6
		\end{pmatrix}&
		\begin{pmatrix}
			6 \\ -2
		\end{pmatrix}&
		\begin{pmatrix}
			-2 \\ 4
		\end{pmatrix}		
	\end{pmatrix}.
\]
The optimal strategies computed by our algorithm with stepsize $t=1/250$ are shown in Figure \ref{fig:01}. Using stepsize $t=1/10$ for I and $t=1/5$ for II we obtain $7$ minimal strategies for I and 5 maximal strategies for II, compare Figures \ref{fig:01} and \ref{fig:02}. Among the resulting 35 pairs there are 10 set Shapley equilibria, two of them are tested numerically to be strong set Shapley equilibria, see Table \ref{tab:t1}. 
In Figure \ref{fig:02}, the relation to POSS payoffs is shown.

\begin{table}[hpt]
	\begin{center}
		\begin{tabular}{lll}
		\toprule 		   
		  $p^T$ & $q^T$ & type \\  
		\midrule		
			$\of{\frac{2}{5}, 0, \frac{3}{5}}$	& $(0,0,1)$ & strong \\[.1cm]
		 	$\of{\frac{1}{2},0,\frac{1}{2}}$	& $(0,0,1)$ & strong \\[.1cm]
			$\of{\frac{3}{5}, 0, \frac{2}{5}}$& $(0,0,1)$ & not strong \\[.1cm]
			$\of{\frac{7}{10}, 0, \frac{2}{10}}$& $(0,0,1)$ & not strong \\[.1cm]
			$\of{\frac{1}{2}, 0, \frac{1}{2}}$	& $\of{\frac{1}{5}, 0, \frac{4}{5}}$& not strong \\
		\bottomrule
		\end{tabular}\hspace{.4cm}
		\begin{tabular}{lll}
		\toprule 		   
		  $p^T$ & $q^T$ & type \\  
		\midrule		
			$\of{\frac{3}{5}, 0, \frac{2}{5}}$	 & $\of{\frac{1}{5}, 0, \frac{4}{5}}$ & not strong \\[.1cm]
			$\of{\frac{7}{10}, 0, \frac{2}{10}}$ & $\of{\frac{1}{5}, 0, \frac{4}{5}}$ & not strong \\[.1cm]
			$\of{\frac{1}{2}, 0, \frac{1}{2}}$	 & $\of{\frac{2}{5}, 0, \frac{3}{5}}$ & not strong \\[.1cm]
			$\of{\frac{3}{5}, 0, \frac{2}{5}}$	 & $\of{\frac{2}{5}, 0, \frac{3}{5}}$  & not strong \\[.1cm]
			$\of{\frac{7}{10}, 0, \frac{2}{10}}$ & $\of{\frac{2}{5}, 0, \frac{3}{5}}$  & not strong \\
		\bottomrule
		\end{tabular}
		\caption{Set Shapley equilibrium points for Example \ref{ex:01} for stepsize $t=1/10$ for I and $t=1/5$ for II.}
		\label{tab:t1}
	\end{center}
\end{table}

\begin{figure}[tp]
	\begin{center}
		\setlength{\unitlength}{0.01\textwidth}%
		\begin{picture}(100,45)
		\includegraphics[width=0.4\textwidth]{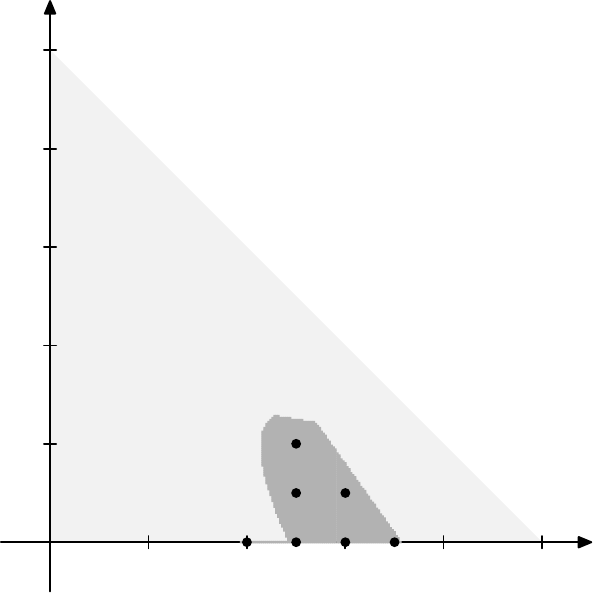}\hspace{.15\textwidth}
		\includegraphics[width=0.4\textwidth]{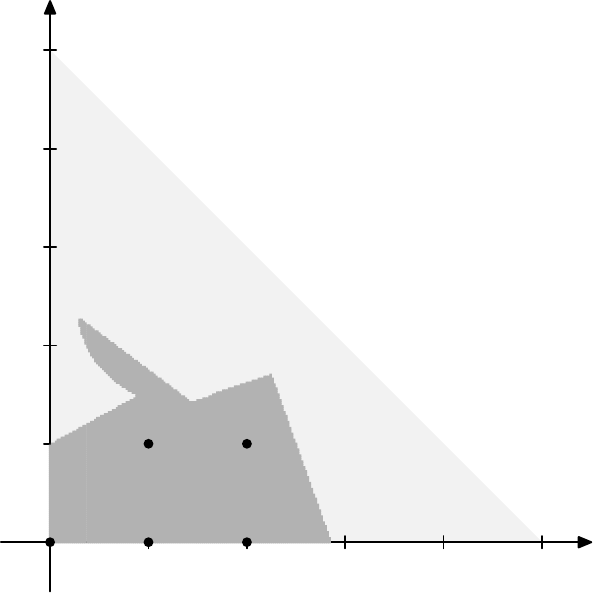}
		\put(-4,0){\small$1$}
		\put(-60,0){\small$1$}
		\put(-39,36){\small$1$}
		\put(-95,36){\small$1$}
		\put(-1,5){$q_1$}
		\put(-56,5){$p_1$}
		\put(-36,41){$q_2$}
		\put(-91.5,41){$p_2$}
		 \end{picture}
		\end{center}
	\caption{Optimal strategies of Example \ref{ex:01} for I (left) and II (right). The first two components $p_1,p_2$ and $q_1,q_2$, respectively, of the strategies $p \in P \subseteq \R^3$ and $q \in Q \subseteq \R^3$ are depicted. The dark gray area refers to optimal strategies and the light gray to non-optimal. The black points are selected strategies the payoff of which is depicted in Figure \ref{fig:02}.}
	\label{fig:01}
\end{figure}
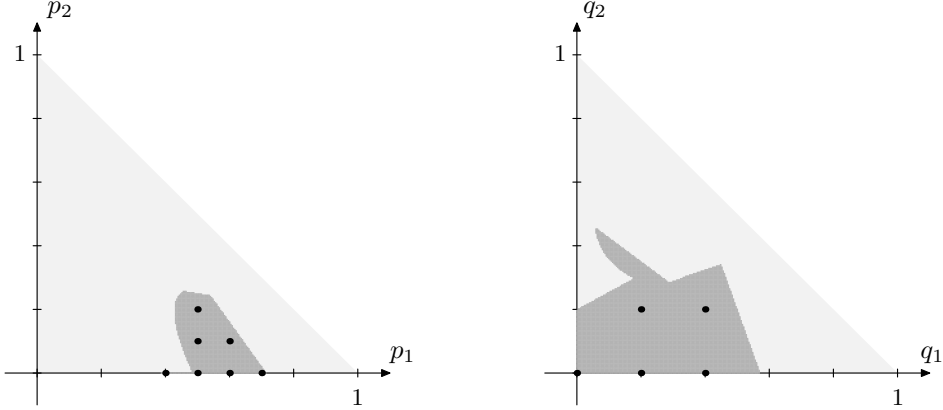

\begin{figure}[t]
	\setlength{\unitlength}{0.01\textwidth}%
	\begin{picture}(100,45)
		\includegraphics[width=0.4\textwidth]{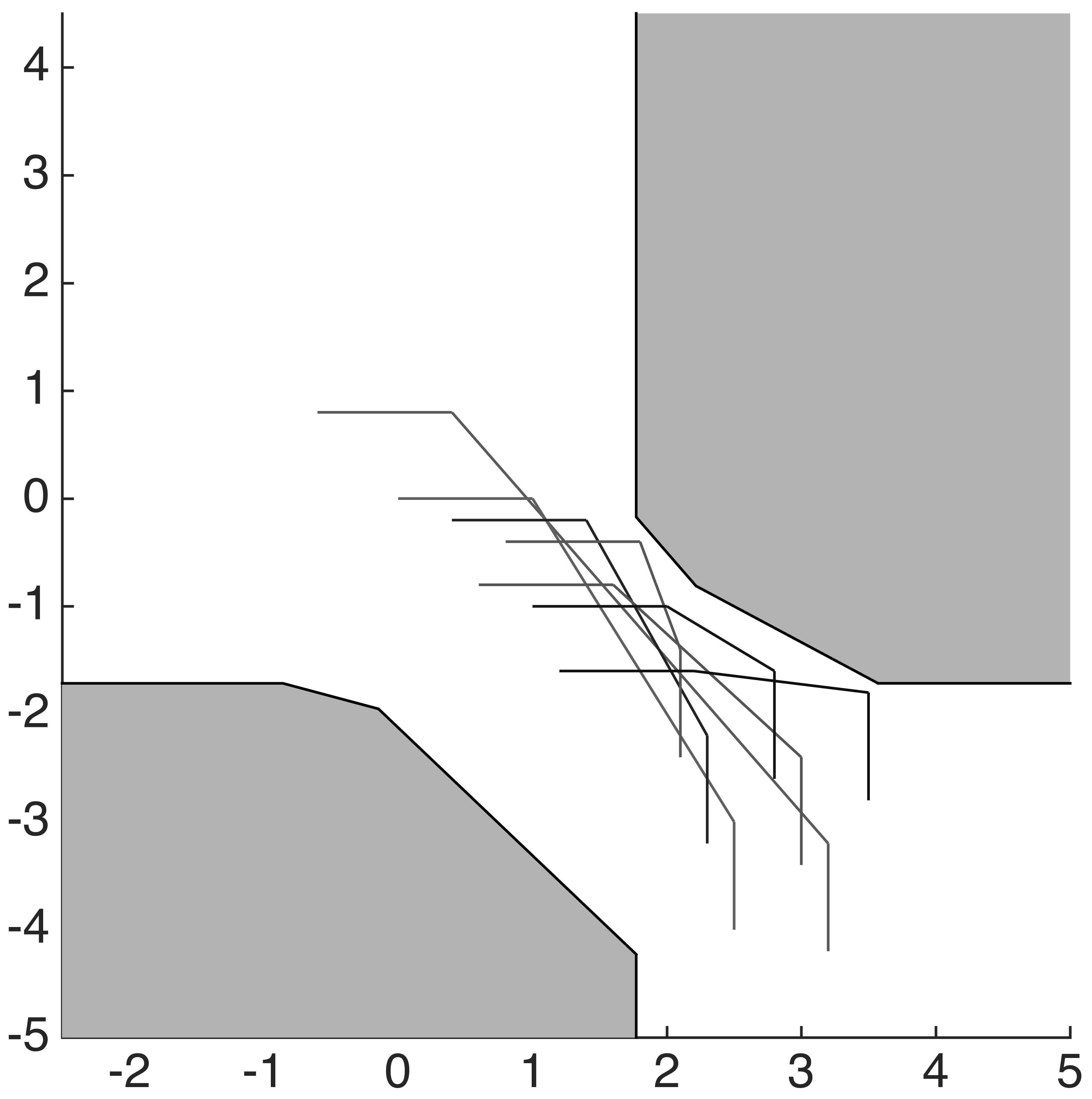}\hspace{.15\textwidth}
		\includegraphics[width=0.4\textwidth]{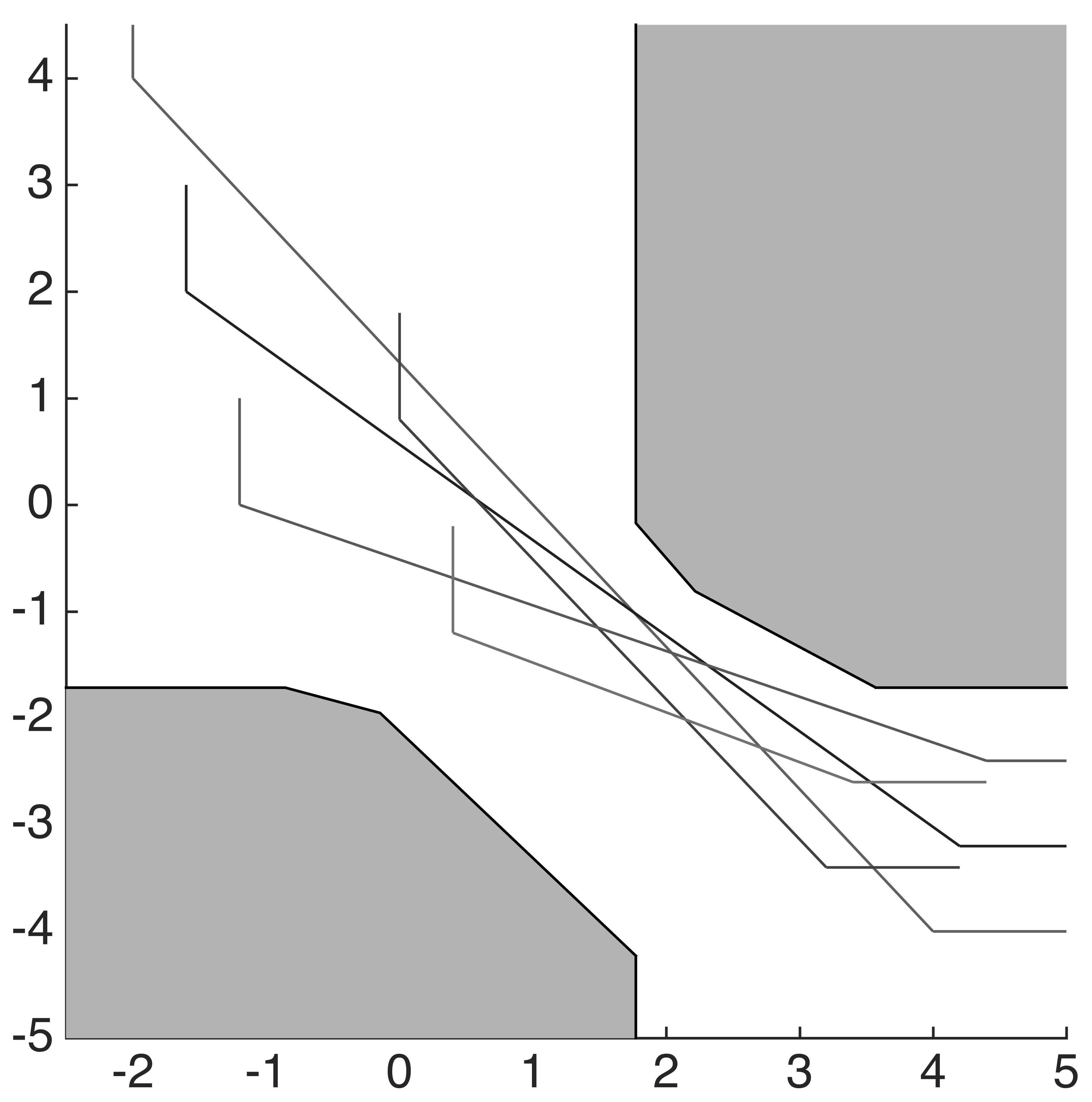}
		\put(-11,30){$\W_I$}
		\put(-69,30){$\W_I$}
		\put(-33,8){$\W_{II}$}
		\put(-91,8){$\W_{II}$}
	\end{picture}
	\begin{picture}(100,45)
		\includegraphics[width=0.4\textwidth]{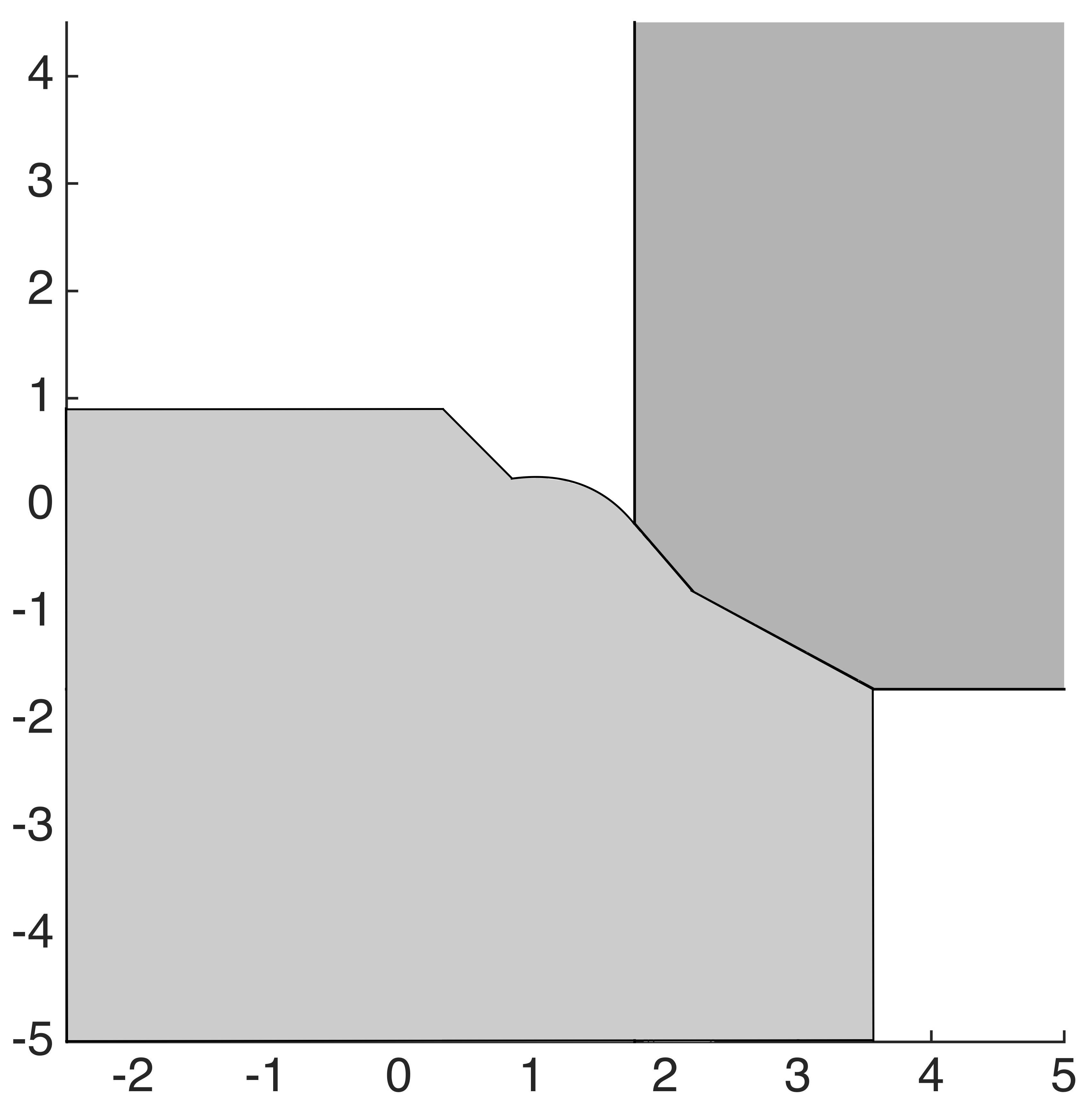}\hspace{.15\textwidth}
		\includegraphics[width=0.4\textwidth]{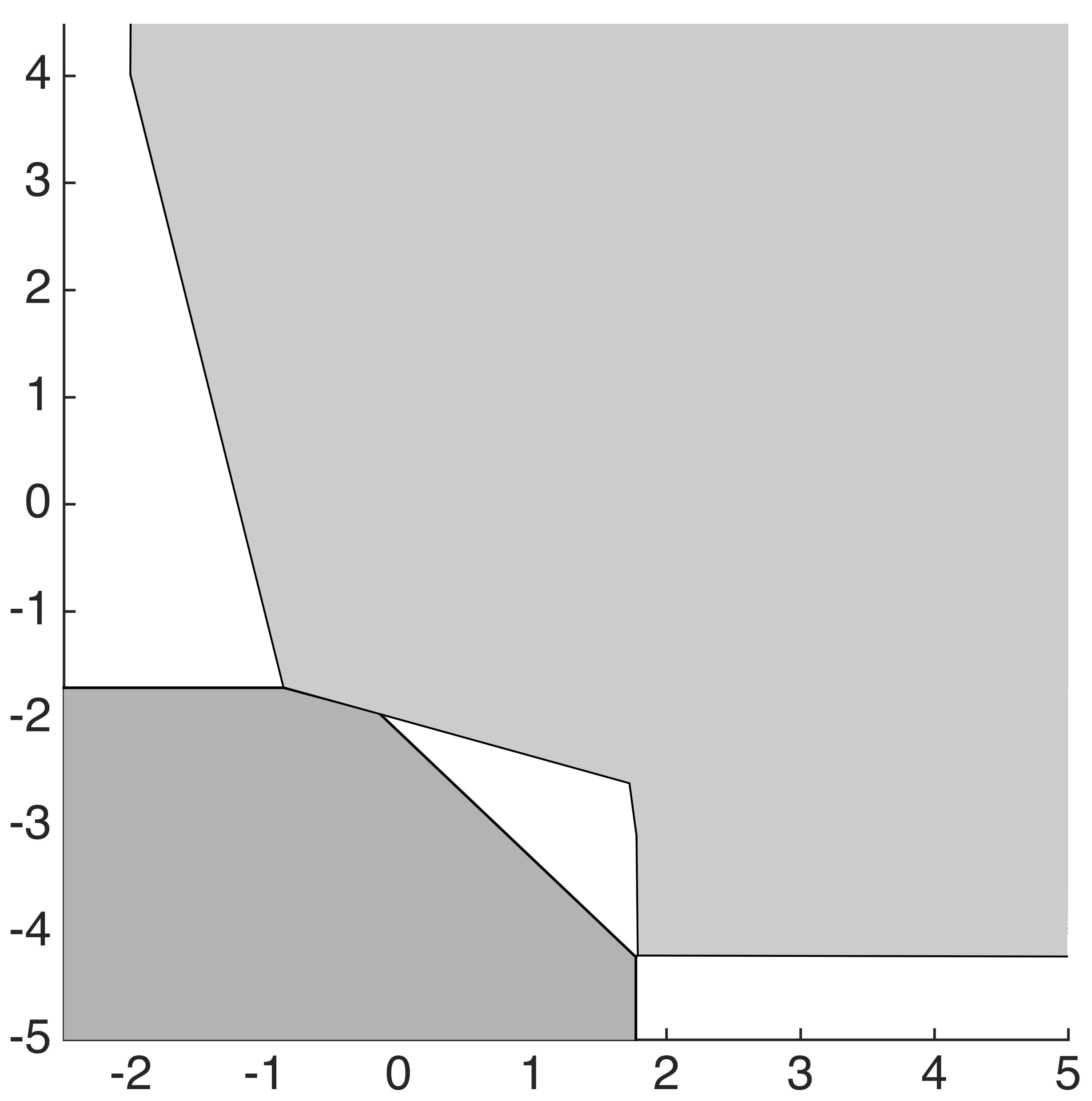}
		\put(-29,28){$\displaystyle\bigcup_{q \in \MAX} V_{II}(q)$}
		\put(-34,8){$\W_{II}$}
		\put(-90,12){$\displaystyle\bigcup_{p \in \MIN} V_{I}(p)$}
		\put(-67,30){$\W_I$}
		\end{picture}
\caption{Payoffs for Example \ref{ex:01}. The sets $\W_I$ and $\W_{II}$ are the payoffs corresponding to POSS strategies. Top left, a selection of sets $V_I(p)$ for optimal strategies $p \in \MIN$ is depicted (only their boundaries are partially shown). The optimal strategies refer to the seven black points in Figure \ref{fig:01} (left). Top right, the same is shown for player II, and the picture is related to Figure \ref{fig:01} (right). By drawing the payoffs for all optional strategies (bottom), one can see that, as stated in Theorem \ref{ThmMinMaxPOSS}, the new approach leads to better results than POSS as even a ``gap'' between the sets can be observed.}
\label{fig:02}
\end{figure}
\end{example}

\begin{example}\label{ex:02}		
Consider a game where all components of $g_{ij} \in \R^K$ are random integer payoffs between -10 and 10. We consider various such games with $m \in\cb{2,3,4,5}$ rows, $n \in \cb{2,3,4,5}$ columns and solve these games by the algorithm described in the previous section. For each choice of $m,n,K$ we solve $5$ random instances. We use a discretization stepsize of $t=1/20$. The average running times for computing optimal strategies (without checking if any pairs are set equilibrium points) are displayed in Table \ref{tab:01}. We observed that different instances of the same dimension may lead to a large variety of percentage of optimal strategies among all strategies considered. However, this had no significant influence on the running times to compute optimal strategies. In the larger examples ($K,m,n \geq 4$) the number of optimal strategies often exceeds $10^3$, hence more than $10^6$ (sometimes even more than $8\cdot 10^7$) pairs have to be checked for being set Shapley equilibria and strong set Shapley equilibria. For $10^5$ pairs of our largest random example ($K=5, m=5, n=5$) we had a running time of $158$ seconds.  
	
\begin{table}[htp]
	\begin{center}
		\begin{tabular}{r@{\hskip 1.2cm}rrrr}
		\toprule 
		  K=2  & \multicolumn{4}{c}{$n$} \\ \cmidrule(r){2-5} 
		  $m$ & 2& 3& 4 & 5\\  
		\midrule		
			2	&	0.2		&		&			&			\\
		 	3	&	0.6		& 1.7	&			&			\\
		 	4	&	3.9		& 7.8	&	14 	&			\\
			5	&	18.8	& 23.5	&	43 	&	78 	\\
		\bottomrule		
		\toprule 
	    K=4 & \multicolumn{4}{c}{$n$} \\ \cmidrule(r){2-5} 
	  $m$ & 2& 3& 4 & 5\\
		\midrule    
		2	&	0.4		&			&			&			\\
 		3	&	2.0		& 	4.7		&			&			\\
		4	&	14.4	&	20.2	&	45	&			\\
		5	&	68.1	& 	105.8	&	186	&	358	\\
		\bottomrule
		\end{tabular}\hfill	
		\begin{tabular}{r@{\hskip 1.2cm}rrrr}
			\toprule 
		    K=3 & \multicolumn{4}{c}{$n$} \\ \cmidrule(r){2-5} 
		  $m$ & 2& 3& 4 & 5\\
			\midrule    
			2	&	0.3		&			&		 	&			\\
	 		3	&	1.6 	& 	3.4		&		 	&			\\
			4	&	7.6		&	10.0	&	30 	&			\\
			5	&	47.1	& 	62.0	&	102 	&	186 	\\
			\bottomrule		
  			\toprule 
  		    K=5& \multicolumn{4}{c}{$n$} \\ \cmidrule(r){2-5} 
  		  $m$ & 2& 3& 4 & 5\\
  			\midrule    
  			2	&	0.4		&			&			&			\\
  	 		3	&	2.0		& 	5.7		&			&			\\
  			4	&	16.1	&	29.6	&	65	&			\\
  			5	&	99.6	& 	136.5	&	238	&	514	\\
  			\bottomrule
		\end{tabular}
		\caption{Running time (in seconds) for Example \ref{ex:02}. The discretization stepsize is $t=1/20$. The average time of 5 random examples is displayed.}
		\label{tab:01}
	\end{center}
\end{table}
\end{example}	

\section{Conclusions and perspectives}

A new solution concept for zero-sum matrix games has been introduced which transfers the notions of minimax- and maximin-strategies from the one-dimensional to the multi-dimensional payoff case; it yields interchangeable strategy pairs and worst case estimates which should be played if the players are ``loss averse" and do not know anything about their preferences, but the fact that they prefer ``less loss" and ``more gain." The new concept uses set relations, but is based on the complete-lattice approach to set optimization. Combining our minimal/maximal solutions with (strengthened) versions of Nash-type equilibrium concepts for multi-dimensional payoff games introduced by Shapley, we obtain new equilibrium concepts and show existence.

Extensions are now possible to situations in which both players have preferences expressed by two potentially different convex cones $C_I, C_{II} \subset \R^K$. In fact, this seems to be just a mathematical exercise since the corresponding concepts, in particular set relations generated by arbitrary cones, are available (see Hamel et al \cite{HamelEtAl15Incoll}). Moreover, even the general situation as considered in Bade \cite{Bade05ET} is well within
reach since every preorder can be extended to set relations (not just vector preorders).

\bibliographystyle{plain}


\begin{thebibliography}{9999}

\bibitem{Aumann62Economet} Aumann, R.J., Utility theory without the completeness axiom, Econometrica 30(3):445-462, 1962

\bibitem{Bade05ET} Bade, S., Nash equilibrium in games with incomplete preferences, Economic Theory 26(2):309-332, 2005
 
\bibitem{Bewley86R} Bewley, T.F., Knightian decision theory. Part I, Discussion Paper no. 807 of the Cowles Foundation at Yale University, 1986
 
\bibitem{Blackwell56PJM} Blackwell, D., An analog of the minimax theorem for vector payoffs, Pacific Journal of Mathematics 6(1):1-8, 1956

\bibitem{BosiHerden12JME} Bosi, G. and Herden, G., Continuous multi-utility representations of preorders, Journal of Mathematical Economics 48(4):212-218, 2012

\bibitem{CampiOwen11FS} Campi, L. and Owen, M.P., Multivariate utility maximization with proportional transaction costs, Finance and Stochastics 15(3):461-499, 2011

\bibitem{CarlierDana13JET} Carlier, G. and Dana, R.-A., Pareto optima and equilibria when preferences are incompletely known, Journal of Economic Theory 148(4): 1606-1623

\bibitem{Cook76NRLQ} Cook, W.D., Zero-sum games with multiple goals, Naval Research Logistics Quarterly 23(4):615-621, 1976

\bibitem{Corley85JOTA} Corley, H.W., Games with vector payoffs, Journal of Optimization Theory and Applications 47(4):491-498, 1985

\bibitem{DeMarcoMorgan07IGTR} De Marco, G. and Morgan, J., A refinement concept for equilibria in multicriteria games via stable scalarizations, International Game Theory Review 9(2):169-181, 2007

\bibitem{DubraMaccheroniOk04JET} Dubra, J. and Maccheroni, F. and Ok, E.A., Expected utility theory without the completeness axiom, Journal of Economic Theory 115(1):118-133, 2004

\bibitem{EliazOk06GEB} Eliaz, K. and Ok, E.A., Indifference or indecisiveness? Choice-theoretic foundations of incomplete preferences, Games and Economic Behavior 56(1):61-86, 2006

\bibitem{Ehrgott05} Ehrgott, M., {\em Multicriteria Optimization}, 2nd edition, Springer-Verlag Berlin 2005

\bibitem{EvrenOk11JME} Evren, {\"O.} and Ok, E.A., On the multi-utility representation of preference relations, Journal of Mathematical Economics 47(4):554-563, 2011

\bibitem{FernandezPuerto96JOTA} Fern{\'a}ndez, F.R. and Puerto, J., Vector linear programming in zero-sum multicriteria matrix games, Journal of Optimization Theory and Applications 89(1):115-127, 1996

\bibitem{FernandezMonroyPuerto98JOTA} Fern{\'a}ndez, F.R. and Monroy, L. and Puerto, J., Multicriteria goal games, Journal of Optimization Theory and Applications 99(2):403-421, 1998

\bibitem{Ghose91JOTA} Ghose, D., A necessary and sufficient condition for Pareto-optimal security strategies in multicriteria matrix games, Journal of Optimization Theory and Applications 68(3):463-481, 1991

\bibitem{GhosePrasad89JOTA} Ghose, D. and Prasad, U.R., Solution concepts in two-person multicriteria games, Journal of Optimization Theory and Applications 63(2):167-189, 1989

\bibitem{GreenMasCollelWhinston95Book} Mas-Colell, A. and Whinston, M.D. and Green, J.R., {\em Microeconomic Theory}, Oxford University Press New York 1995



\bibitem{HamelEtAl15Incoll} Hamel, A.H. and Heyde, F. and L{\"o}hne, A. and Rudloff, B. and Schrage, C., Set optimization--a rather short introduction. In: Hamel, A.H. and Heyde, F. and L{\"o}hne, A. and Rudloff, B. and Schrage, C. (eds.), {\em Set optimization and applications -- the state of the art. From set relations to set-valued risk measures}, Springer-Verlag Berlin 2015, pp. 65-141

\bibitem{Henig86} Henig, M.I., The domination property in multicriteria optimization, Journal of Mathematical Analysis and Applications 114(1):7-16, 1986

\bibitem{HeydeLoehne11Opt} Heyde, F. and L\"ohne, A., Solution concepts in vector optimization: a fresh look at an old story, Optimization 60(10-12):1421-1440, 2011

\bibitem{KuroiwaTanakaTruong97NA} Kuroiwa, D. and Tanaka, T. and Ha, T.X.D., On cone convexity of set-valued maps, Nonlinear Analysis. Theory, Methods \& Applications 30(3):1487-1496, 1997

\bibitem{LoeWei16} L{\"o}hne, A. and Wei{\ss}ing, B., The vector linear program solver Bensolve -- notes on theoretical background, European Journal of Operational Research 2016, DOI: 10.1016/j.ejor.2016.02.039

\bibitem{LoeWei15} L{\"o}hne, A. and Wei{\ss}ing, B., Equivalence between polyhedral projection, multiple objective linear programming and vector linear programming, Mathematical Methods of Operations Research 84(2):411-426, 2016
	
\bibitem{bensolve} L\"ohne, A. and Wei{\ss}ing, B., Bensolve - {VLP} solver, version 2.0.1, www.bensolve.org, unpublished

\bibitem{LucVargas92NA} Luc, D.T. and Vargas, C., A saddlepoint theorem for set-valued maps, Nonlinear Analysis: Theory, Methods \& Applications 18(1):1-7, 1992

\bibitem{Maeda15Incoll} Maeda, T., On characterization of Nash equilibrium strategy in bi-matrix games with set payoffs. In: Hamel, A.H. and Heyde, F. and L{\"o}hne, A. and Rudloff, B. and Schrage, C. (eds.), {\em Set optimization and applications -- the state of the art. From set relations to set-valued risk measures}, Springer-Verlag Berlin 2015, pp. 313-331

\bibitem{Nieuwenhuis83JOTA} Nieuwenhuis, J.W., Some minimax theorems in vector-valued functions, Journal of Optimization Theory and Applications 40(3):463-475, 1983

\bibitem{Ok02JET} Ok, E.A., Utility representation of an incomplete preference relation, Journal of Economic Theory 104(2):429-449, 2002

\bibitem{OkOrtolevaRiella12Economet} Ok, E.A. and Ortoleva, P. and Riella, G., Incomplete preferences under uncertainty: Indecisiveness in beliefs versus tastes, Econometrica 80(4):1791-1808, 2012

\bibitem{Shapley59NRLQ} Shapley, L.S., Equilibrium points in games with vector payoffs, Naval Research Logistics Quarterly 6:57-61, 1959
 
\bibitem{Tanaka88JOTA} Tanaka, T., Some minimax problems of vector-valued functions, Journal of Optimization Theory and Applications 59(3):505-524, 1988

\bibitem{Tanaka94JOTA} Tanaka, T., Generalized quasiconvexities, cone saddle points, and minimax theorem for vector-valued functions, Journal of Optimization Theory and Applications 81(2):355-377, 1994

\bibitem{Tanaka00Incoll} Tanaka, T., Vector-Valued Minimax Theorems in Multicriteria Games. In: Yong Shi and Zeleny, M. (eds.), {\em New Frontiers of Decision Making for the Information Technology Era}, World Scientific 2000, pp. 75-99

\bibitem{Wierzbicki95JSEE} Wierzbicki, A.P., Multiple criteria games--theory and applications, Journal of Systems Engineering and Electronics 6(2):65-81, 1995

\bibitem{Zeleny75IJGT} Zeleny, M., Games with multiple payoffs, International Journal of Game Theory 4(4):179-191, 1975

\bibitem{Zeleny74} Zeleny, M., {\em Linear Multiobjective Programming}, Lecture Notes in Economics and Mathematical Systems, Vol. 95, Springer-Verlag, Berlin-New York, 1974

\bibitem{Zhao91IJGT} Zhao, J., The equilibria of a multiple objective game, International Journal of Game Theory 20(2):171-182, 1991

\end{thebibliography}

\section*{Appendix}

For the readers convenience, the appendix summarizes basic concepts related to vector and set orders as well as the complete lattice approach to set optimization.

Let $C \subseteq \R^K$ be a closed convex cone satisfying $C \cap (-C) = \cb{0}$. Such a cone generates a partial order on $\R^K$ (i.e. a reflexive, transitive and antisymmetric relation) by
\[
y \leq_C z \quad \Leftrightarrow \quad z - y \in C,
\]
and this order is compatible with the linear space operations on $\R^K$, i.e. $\leq_C$ is a vector order. In general, $(\R^K, \leq_C)$ is not a lattice. Even if the infimum (or the supremum) of a set $A \subseteq \R^K$ exists, it can be ``far away" from $A$: consider $C = \R^2_+$ and $A =\cb{z \in \R^2_+ \mid z_1+z_2 \geq 2}$ whose infimum with respect to $\leq_{\R^2_+}$ is $z = 0 \in \R^2$. Therefore, the predominant optimality notion in vector optimization and multi-criteria decision making is based on minimal (or maximal) points.

A point $\bar z \in A \subseteq \R^K$ is called minimal with respect to $\leq_C$ if
\[
z \in A, \; z \leq_C \bar z \quad \Rightarrow \quad z = \bar z.
\]
The set of minimal points of $A$ is denoted by $\Min A$. Likewise, the set $\Max A$ of maximal points is introduced.

The lack of a reasonable infimum/supremum with respect to vector orders is a major motivation for introducing so-called set relations, see Kuroiwa et al \cite{KuroiwaTanakaTruong97NA} as well as Hamel et al \cite{HamelEtAl15Incoll} for a recent survey with many references.

Let $A, B \subseteq \R^K$. By
\[
A \lel_C B  \; :\Leftrightarrow \; B \subseteq A + C \quad \text{and} \quad
A \leu_C B  \; :\Leftrightarrow \; A \subseteq B - C
\]
two ``set relations" are defined which both are reflexive and transitive, but not antisymmetric in general. Moreover, they are two different extensions of $\leq_C$: $z \leq_C y$ $\Leftrightarrow$ $\cb{z} \lel_C \cb{y}$  $\Leftrightarrow$ $\cb{z} \leu_C \cb{y}$. Per se, these relations just shift the difficulty of defining optimality from $\R^K$ to its power set $\mathcal P(\R^K)$. Their value lies in the possibility to construct complete lattices of sets based on their symmetric parts.

Two sets $A, B \subseteq \R^K$ are equivalent with respect to $\lel_C$, written $A \sim B$, if $A \lel_C B \lel_C A$. It can easily be shown that $A \sim B$ if and only if $A + C = B + C$. Therefore, the set of equivalence classes with respect to $\sim$ can be identified with $\mathcal P(\R^K, C) := \cb{A \subseteq \R^K \mid A = A + C}$. Moreover, on $\mathcal P(\R^K, C)$ the relation $\lel_C$ coincides with $\supseteq$. Likewise, $\mathcal P(\R^K, -C) := \cb{A \subseteq \R^K \mid A = A - C}$ can be identified with the set of equivalence classes with respect to the symmetric part of $\leu_C$, and on $\mathcal P(\R^K, -C)$ the relation $\leu_C$ coincides with $\subseteq$.

Moreover, both $\of{\mathcal P(\R^K, C), \supseteq}$ and $\of{\mathcal P(\R^K, -C), \subseteq}$ are complete lattices, i.e. every set $\mathcal A \subseteq \mathcal P(\R^K, C)$ has an infimum and a supremum in $\mathcal P(\R^K, C)$ as well as every set $\mathcal B \subseteq \mathcal P(\R^K, -C)$. Thus, infimum and supremum become available again without any restrictions, and this fact constitutes the major difference to more traditional approaches in multi-objective and even set-valued optimization, the latter only based on minimality/maximality with respect to the set relations defined above. Compare the survey Hamel et al \cite{HamelEtAl15Incoll} for more details and references.

The applications in this paper involve convex (set-valued) functions and closed convex sets. Therefore, the following two sets are introduced:
\begin{align*}
\G(\R^K, C) & := \cb{A \subseteq \R^K \mid A = \cl\co(A + C)} \\
\G(\R^K, -C) & := \cb{A \subseteq \R^K \mid A = \cl\co(A - C)}
\end{align*}
where the usual conventions for the Minkowski addition of sets are used with the extension $\emptyset + A = A + \emptyset$ for all $A \in \P(\R^K)$.

\begin{proposition}
\label{PropComLatt}
The pairs $\of{\G(\R^K, C), \supseteq}$ and $\of{\G(\R^K, -C), \subseteq}$ are complete lattices with the following formulas for infimum and supremum: For $\mathcal A \subseteq \G(\R^K, C)$,
\[
\inf \mathcal A = \cl\co\bigcup_{A \in \mathcal A} A \quad \text{and} \quad
	\sup \mathcal A = \bigcap_{A \in \mathcal A} A.
\]
For $\mathcal B \subseteq \G(\R^K, -C)$,
\[
\inf \mathcal B = \bigcap_{B \in \mathcal B} B \quad \text{and} \quad
	\sup \mathcal B = \cl\co\bigcup_{B \in \mathcal B} B.
\]
\end{proposition}

{\sc Proof.} See, for example, \cite{HamelEtAl15Incoll}. \pend

\medskip Note that the formulas for inf and sup in $(\mathcal G\of{\R^K, C}, \supseteq)$ and $(\mathcal G\of{\R^K, -C}, \subseteq)$ are exchanged due to the change of the ordering relation.


A set $A \subseteq \R^K$ is said to enjoy the lower domination property if for each $z \in A$ there is $\bar z \in \Min A$ with $\bar z \leq_C z$. The upper domination property is defined parallel.

If $A, B \subseteq \R^K$ satisfy the lower domination property, then
\[
A \lel_C B \quad \Leftrightarrow \quad \Min B \subseteq \Min A + C,
\] 
and if they satisfy the upper domination property, then
\[
A \leu_C B \quad \Leftrightarrow \quad \Max A \subseteq \Max B - C.
\]
It is well-known that $A \subseteq \R^K$ satisfies the lower as well as the upper domination property if it is compact, see e.g. \cite{Henig86}.


The complete-lattice approach admits to provide precise solution concepts for optimization problems with a vector- or set-valued objective function. 

First, we provide a few basic concepts and facts on functions mapping into complete lattices of sets. A function $f \colon \R^k \to \G(\R^K, C)$ is called convex if $s \in (0,1)$ and $x,y \in \R^k$ imply
\[
f(sx + (1-s)y) \supseteq sf(x) + (1-s)f(y).
\]
If $f \colon \R^k \to \G(\R^K, C)$ is convex, then
\[
\inf_{x \in \R^k} f(x) = \cl\bigcup_{x \in \R^k}f(x)
\]
since the set $\bigcup_{x \in \R^k}f(x)$ already is convex as one easily checks.

A function $g \colon \R^k \to \G(\R^K, -C)$ is called concave if $s \in (0,1)$ and $x,y \in \R^k$ imply
\[
sg(x) + (1-s)g(y) \subseteq g(sx + (1-s)y).
\]
Again, if $g$ is concave, the set $\bigcup_{x \in \R^k}g(x)$ is convex and the convex hull in the formula for the $\sup_{x \in \R^k} g(x)$ can be dropped.

Finally, a solution concept for set optimization problems is given which is due to Heyde, L\"ohne \cite{HeydeLoehne11Opt}.

\begin{definition}
\label{DefLatticeSolution}
Let $X$ be a nonempty set, $(L, \leq)$ a complete lattice and $f \colon X \to L$ a function.

(a) A set $M \subseteq X$ is called an infimizer for $f$ if
\[
\inf_{x \in M} f(x) = \inf_{x \in X} f(x).
\]

(b) A point $\bar x \in X$ is called a minimizer for $f$ if
\[
x \in X, \; f(x) \leq f(\bar x) \quad \Rightarrow \quad f(x) = f(\bar x).
\]

(c) A set $M \subseteq X$ is called a solution of the problem
\[
\tag{P} \text{minimize} \quad f(x) \quad \text{over} \quad x \in X
\]
if $M$ is an infimizer and each $x \in M$ is a minimizer for $f$. A solution $M$ is called full, if $M$ includes all minimizers for $f$. 
\end{definition}

\end{document}